\documentclass[11pt, oneside, reqno]{amsart}
\usepackage{amsfonts, amssymb, color}
\usepackage{amsmath}
\usepackage{graphicx}
\usepackage{amscd}
\usepackage{amsthm}
\usepackage{tikz}
\usepackage{subcaption}
\usepackage{hyperref}
\usepackage{longtable}

\usetikzlibrary{decorations.markings}
\usetikzlibrary{shapes,arrows,positioning}	
\usetikzlibrary {arrows.meta}	
\usetikzlibrary{decorations.markings,decorations.pathmorphing}
\usetikzlibrary{decorations.pathreplacing}
\usetikzlibrary{calc, fit}

\newcommand{\rmidarrow}{\tikz \draw[thick,-{Computer Modern Rightarrow}] (-1pt,0) --(1pt,0);}	
\newcommand{\lmidarrow}{\tikz \draw[thick, -{Computer Modern Rightarrow}] (1pt,0) -- (-1pt,0);}

\newcommand{\Nwmidarrow}{\tikz \draw[thick, -{Computer Modern Rightarrow}] (0,0) -- (2pt,2pt);}
\newcommand{\Nmidarrow}{\tikz \draw[thick, -{Computer Modern Rightarrow}] (0,0) -- (0,1);}
\newcommand{\Nemidarrow}{\tikz \draw[thick, -{Computer Modern Rightarrow}] (0,0) -- (-2pt, 2pt);}
\usepackage[margin=1in]{geometry}

\newcommand{\mystar}{\operatornamewithlimits{\ast}}

	\newtheorem{thm}{Theorem}[section]
	
	\newtheorem{lem}[thm]{Lemma}
	\newtheorem{prop}[thm]{Proposition}
	\newtheorem{cor}[thm]{Corollary}
	
	\newtheorem{rem}[thm]{Remark}

	\newtheorem{ex}[thm]{Example}

\begin{document}
	
\title[]{Differences Between Quantum Symmetries of Graphs and Graph $C^*$-algebras}

\author{Ujjal Karmakar}
\address{Stat-Math Unit, Indian Statistical Institute, 203, B. T. Road, Kolkata-700108, India}
\email{mathsujjal@gmail.com}

\author{Arnab Mandal}
\address{Department of Mathematics, Presidency University, College Street, Kolkata-700073, India}
\email{arnab.maths@presiuniv.ac.in}

\maketitle

\begin{abstract}
We establish several characteristic differences between the quantum symmetry groups of directed graphs and their associated graph $C^*$-algebras in the frameworks of Banica and Joardar-Mandal, respectively, considering the family of finite directed graphs without isolated vertices and without multiple edges.
\end{abstract}
 

\section{{\bf Introduction}}
In 2003, Bichon first introduced the notion of the quantum symmetry group $QAut_{Bic}(\Gamma)$ of a finite graph $\Gamma$ \cite{Bichon}, generalizing Wang's notion of the quantum symmetry group $S_n^+$ of a finite $n$-point space \cite{Wang}. Later, in 2005, this notion was further developed by Banica \cite{Ban}, obtaining the inclusion $Aut(\Gamma) \subseteq QAut_{Bic}(\Gamma) \subseteq QAut_{Ban}(\Gamma) \subseteq S_n^+$ in the quantum setting. Here, $QAut_{Ban}(\Gamma)$ denotes Banica's version of the quantum symmetry group of a graph $\Gamma$. A rich body of literature on Banica's quantum automorphism group of graphs has already been developed to date (\cite{Banver} -- \cite{BBC}, \cite{Qauttrees} -- \cite{Fulton}, \cite{Lupini}, \cite{solutiongroup}, \cite{SS_petersen} -- \cite{SSthesis} and references therein). \\ 
On the other hand, around 2017-2018, Schmidt and Weber \cite{Web} first introduced the notion of quantum symmetry of graph $C^*$-algebra by observing that the quantum symmetry group of the underlying graph acts maximally on the given graph $C^*$-algebra. Thus, their notion of quantum symmetry of a graph $C^*$-algebra actually coincides with the quantum symmetry group of the underlying graph. Inspired by their work, Joardar and Mandal \cite{Mandal, Mandalkms} came up with another notion of quantum symmetry of the graph $C^*$-algebra (denoted by $Q_{\tau}^{Lin}(\Gamma)$) in a categorical framework in 2018. Interestingly, $Q_{\tau}^{Lin}(\Gamma)$ is always a quantum subgroup of $A_{u^t}(F^{\Gamma})$ for some invertible diagonal matrix $F^{\Gamma}$ and contained the quantum symmetry group $QAut_{Ban}(\Gamma)$ of the underlying graph as a quantum subgroup. This article compares the characteristic behavior of the quantum symmetry groups of a directed graph $\Gamma$ and its associated graph $C^*$-algebras. To make this comparison in a uniform framework, {\bf we consider finite directed graphs without isolated vertices or multiple edges}, so that the notion of two quantum symmetry groups can be defined. Alternatively, one can think of this comparison as both notions of quantum symmetries of a graph $C^*$-algebra in the sense of Schmidt-Weber and Joardar-Mandal.\\ 

One of the important observations is that $\left( \mystar\limits_{|E(\Gamma)|} C(S^1) \right) $ is a quantum subgroup of $Q_{\tau}^{Lin}(\Gamma)$ [Proposition 3.11, \cite{Mandal}], and thus $Q_{\tau}^{Lin}(\Gamma)$ is always non-commutative and infinite dimensional except for two graphs $P_1$ and $L_1$ with exactly one edge. In contrast, given any $n \in \mathbb{N}$, it is possible to find directed graphs $\Gamma$ with $n$ vertices (resp., $n$ edges) whose quantum symmetry group $QAut_{Ban}(\Gamma)$ is commutative and finite-dimensional -- even trivial, namely $\mathbb{C}$. This observation motivates us to identify some further differences between the quantum symmetry groups of graphs and those of their associated graph $C^*$-algebras. Firstly, it can be observed that there is no common quantum group that can appear as the quantum symmetry of a graph as well as the quantum symmetry of a graph $C^*$-algebra [Proposition \ref{no common q. symm.}]. Furthermore, a quantum symmetry group can actually occur as the quantum symmetry group of infinitely many pairwise non-isomorphic connected [Theorem \ref{graph_infinitely many}] as well as disconnected graphs [Corollary \ref{discon. inf. many graphs}], whereas the quantum symmetry group of a graph $C^*$-algebra can arise from only finitely many graphs [Theorem \ref{graph C*alg_finitely many}]. Consequently, we observe that having the same quantum symmetry group $Q_{\tau}^{Lin}(.)$ for two different graphs forces the graphs to have the same number of edges [Corollary \ref{edge equal}]. This reflects a certain rigidity in the behavior of the quantum symmetry groups of graph $C^*$-algebras. However, the quantum symmetry group of a graph does not bound the cardinality of either the vertex set or the edge set (see the construction in Theorem \ref{graph_infinitely many}). It is known that a graph and its complement graph always have the same quantum symmetry; in contrast, the analogous result does not hold in general for the quantum symmetry groups of graph $C^*$-algebras. However, we describe two infinite families of graphs: one for which the isomorphism $Q_{\tau}^{Lin}(\Gamma) \cong Q_{\tau}^{Lin}(\Gamma^c)$ holds [(2) of Section \ref{sec_comp. graph}], and another for which it does not [(1) of Section \ref{sec_comp. graph}]. Next, we focus on how two notions of quantum symmetry groups behave on the disjoint union of two non-isomorphic graphs. It is known that for simple finite undirected graphs, the quantum symmetry group of the disjoint union of two connected quantum non-isomorphic graphs is the same as the free product of the individual quantum symmetry groups, namely $QAut_{Ban}(\Gamma_1 \sqcup \Gamma_2) \cong QAut_{Ban}(\Gamma_1) * QAut_{Ban}(\Gamma_2)$. We first ensure that similar results hold in the setting of finite directed graphs without multiple edges [Proposition \ref{graph_free prod}], and then provide counter-examples of an infinite family of graphs for which analogous isomorphism fails to  hold (even for quantum non-isomorphic graphs) [Proposition \ref{graph C^*alg_free prod not hold}] in the context of quantum symmetry groups of graph $C^*$-algebras (that is, replacing $QAut_{Ban}(.)$ with $Q_{\tau}^{Lin}(.)$). However, a similar isomorphism may also hold in the graph $C^*$-algebraic setting for the family of path graphs $P_n$ with distinct path lengths  [Proposition \ref{matix free prod}]. Finally, we compare the behaviour of the minimal quantum symmetry groups in both aspects. We show that the quantum symmetry of graph $C^*$-algebra is minimal (that is, of the form free product $C(S^1)$) implies the underlying graph is quantum asymmetric [Proposition \ref{rigid imply q.assym}], but the converse is not true.\\

In the second section of this article, we discuss the required preliminaries on the quantum symmetry of graphs and graph $C^*$-algebras and recall the relevant results, retaining the same theorem and proposition numbering as in the original sources. The third section (the main section of this article) is divided into five subsections, each highlighting a characteristic difference between the two notions of quantum symmetry groups. Within a subsection, the discussions regarding the quantum symmetry of graphs and  the quantum symmetry of graph $C^*$-algebras are indicated under the bullet marks `$\blacksquare$' and `$\blacktriangleright$' respectively.

\section{ {\bf Preliminaries}}
\subsection{Notations and conventions}
For $n \in \mathbb{N}$, $[n]$ denotes the set $\{1, 2,..., n\}$.  
The tensor product `$\otimes$' denotes the spatial (or minimal) tensor product between $C^*$-algebras. The matrix $I_n$ denotes the identity matrix of $M_n(\mathbb{C})$. The matrices $0_{m,n}$ and $J_{m,n} \in M_{m \times n}(\mathbb{C})$ denote the matrices with all entries equal to 0 and 1, respectively; in particular, we write $0_{n}:=0_{n,n}$ and $J_{n}:= J_{n,n}$. The notation $G \wr S_n$ denotes the wreath product of the group $G$ with the permutation group $S_n$.\\
All the $C^*$-algebras considered here are assumed to be unital. 

\subsection{Finite graph and its graph $C^*$-algebra}
In this subsection, we recall some basic terminology regarding directed graphs and useful facts about graph $C^{*}$-algebras.\\
 A {\bf finite directed graph} $\Gamma:=\{V(\Gamma), E(\Gamma), s, r \}$  consists of  a finite set $V(\Gamma)$ of vertices and a finite set $E(\Gamma)$ of edges such that together with maps $s,r: E(\Gamma) \to V(\Gamma) $ describe the source and range of the edges. {\bf Throughout this article, we consider only the finite directed graphs, which we will often refer to simply as directed graphs}.  The {\bf order} of a graph $\Gamma$ is defined as the cardinality of its set of vertex $V(\Gamma)$. For $v,w \in V(\Gamma)$, the notation `$v \to w $' means that there exists an edge $e \in E(\Gamma)$ such that  $ v=s(e)$ and $ w=r(e) $. We say that $v$ and $w$ are {\bf adjacent} (denoted by $v \sim w$) if $v \to w$ or $w \to v$. The {\bf indegree} of a vertex $v$ (resp. {\bf outdegree} of a vertex $v$) in $\Gamma$, denoted by $indeg_{\Gamma}(v)$ (resp. $outdeg_{\Gamma}(v)$), is the cardinality of the set $\{e \in E(\Gamma) : r(e)=v\}$ (resp. $\{e \in E(\Gamma) : r(e)=v\}$). If the graph is already specified in the context, we simply denote it by $indeg(v)$ and $outdeg(v)$. A vertex $v$ is called a {\bf sink} if $s^{-1}(v)= \emptyset$. A directed graph \textbf{without isolated vertices} means for every vertex $v \in V(\Gamma) $ at least one of $s^{-1}(v)$ or $r^{-1}(v)$ is nonempty. A \textbf{path} of length $n$ is a sequence $\alpha=e_{1}e_{2} \cdots e_{n} $ of edges in $\Gamma$ such that $ r(e_{i}) = s(e_{i+1}) $ for $1 \leq i \leq (n-1) $, and we define $s(\alpha):=s(e_{1})$ and  $r(\alpha):=r(e_{n})$. The graph $P_n$ denotes a path of length $n$ consisting of $n+1$ vertices. The graph $C_n$ denotes a path of length $n$ consisting of $n$ vertices such that $s(e_1)=r(e_n)$. A \textbf{loop} is an edge $l$ with $s(l)=r(l)$. A graph {\bf without multiple edges} means that there do not exist two distinct edges $e, f \in E(\Gamma)$ such that $s(e)=s(f)$ and $r(e)=r(f)$. 
 A graph is said to be {\bf weakly connected} if any two distinct vertices $v,w \in V(\Gamma)$, either $v \sim w$ or there exists a finite sequence of vertices $v_1, v_2,...,v_n$ such that $v \sim v_1 \sim v_2 \sim \cdots \sim v_n \sim w$. We say that a graph is {\bf disconnected} if it is not weakly connected.
The {\bf adjacency matrix} of a graph $\Gamma=\{V(\Gamma),E(\Gamma),s,r \}$ without multiple edges, with respect to the ordering of the vertices $ (v_{1},v_{2},..., v_{n}) $, is an $n \times n$ matrix $ A_{\Gamma}$ with 
$ (A_{\Gamma})_{ij} =
\begin{cases}
1 & if ~~ v_{i} \to v_{j} \\
0 &  ~~ otherwise 
\end{cases}. $ \\ 

\noindent Given a finite, directed graph $\Gamma$, the {\bf graph $C^{*}$-algebra} $C^{*}(\Gamma)$ is a universal $C^{*}$-algebra generated by partial isometries $ \{S_{e}: e \in E(\Gamma) \} $ and orthogonal projections $ \{p_{v}: v \in V(\Gamma) \}$ such that
	\begin{itemize}
		\item[(CK1)] $S_{e}^{*}S_{e}=p_{r(e)}$ for all $ e \in E(\Gamma) $, \vspace{0.1cm}
		\item[(CK2)] $p_{v}=\sum\limits_{\{f:s(f)=v\}}S_{f}S_{f}^{*}$ for all $ v \in V(\Gamma) $ if $s^{-1}(v) \neq \emptyset $.\\
	\end{itemize}
It is always a unital $C^*$-algebra and $\sum\limits_{v \in V(\Gamma)}p_{v}=1$. We refer the reader to \cite{BPRS}, \cite{MS}, \cite{Pask}, \cite{Raeburn} for details regarding graph $C^*$-algebras.\\

 \subsection{Compact quantum groups }
 In this subsection, some important aspects of compact quantum groups and quantum symmetry groups are briefly recalled, which will be required for this article. We refer the reader to \cite{Van}, \cite{Nesh}, \cite{Tim}, \cite{universalCQG}, \cite{Wang} -- \cite{remark} for further details.\\
 
 A {\bf compact quantum group (CQG)} is a pair $(\mathcal{Q}, \Delta) $, where $\mathcal{Q}$ is a unital $C^{*}$-algebra and the coproduct $ \Delta : \mathcal{Q} \to \mathcal{Q} \otimes \mathcal{Q} $  is a unital $C^{*}$-homomorphism such that
 	\begin{itemize}
 		\item[(i)] $(id_{\mathcal{Q}} \otimes \Delta)\Delta = (\Delta \otimes id_{\mathcal{Q}})\Delta $. 
 		\item[(ii)] $ span\{ \Delta(\mathcal{Q})(1 \otimes \mathcal{Q} )\}$ and $ span\{\Delta(\mathcal{Q})(\mathcal{Q} \otimes 1)\} $ are dense in $(\mathcal{Q}\otimes \mathcal{Q})$.
 	\end{itemize}
Given two CQGs $(\mathcal{Q}, \Delta)$ and $(\mathcal{Q}', \Delta') $, a compact quantum group morphism (CQG morphism) from $(\mathcal{Q}, \Delta)$ to $(\mathcal{Q}', \Delta')$ is a $C^{*}$-homomorphism $ \phi: \mathcal{Q} \to \mathcal{Q}' $ such that $ (\phi \otimes \phi)\Delta=\Delta' \phi $.\\
\noindent The CQG is sometimes denoted by $\mathbb{G}$ and the underlying $C^*$-algebra associated to it by $C(\mathbb{G})$, i.e. $\mathbb{G}=(C(\mathbb{G}),\Delta)$.
 Two CQGs $(A, \Delta)$ and $(A', \Delta')$ are said to be {\bf CQG-isomorphic} (simply denoted by $A \cong A'$) if there exists a $C^{*}$-isomorphism $ \phi: A \to A' $ such that $ (\phi \otimes \phi)\Delta=\Delta' \phi $. More generally, a CQG $(A',\Delta')$ is said to be a {\bf quantum subgroup} of $(A, \Delta)$ if there exists a surjective $C^*$-homomorphism $ \phi: A \to A' $ such that $ (\phi \otimes \phi)\Delta=\Delta' \phi $. \\
 An important class of CQG is known as {\bf compact matrix quantum group (CMQG)}, denoted by $(A,q)$, which is by definition a unital $C^*$-algebra $A$ together with a matrix $q=(q_{ij})_{n \times n}$ whose entries generate $A$ and satisfy the following:
 	\begin{enumerate}
 		\item[(i)] both $q$ and ${q}^t:=(q_{ji})_{n \times n}$ are invertible in $M_n(A)$,
 		\item[(ii)] there exists a unital $C^*$-homomorphism $\Delta: A \to A \otimes A $, referred to as the coproduct, such that $\Delta(q_{ij})= \sum_{k=1}^{n}q_{ik} \otimes q_{kj}$ for all $i,j \in \{1,2,...,n\}$. 
 	\end{enumerate}
 	The matrix $q$ is called the {\bf fundamental representation} of $A$. 
 In a CMQG $(A,q)$, if $A_0$ denotes a *-subalgebra generated by $\{q_{ij}\}_{i,j \in [n]}$, then there exists a linear anti-multiplicative map $\kappa: A_0 \to A_0$, called the {\bf antipode}, such that $\kappa(q_{ij})=(q^{-1})_{ij}$, where $q^{-1}$ is the inverse of $q$. \\

\noindent The following theorem ensures the existence of a CQG morphism that factors through a quotient of a CQG by an ideal.\\
\noindent {\bf [Theorem 2.11, \cite{Wangfree}]:} For a Woronowicz $C^*$-ideal $J$ of $(\mathcal{Q}, \Delta)$, the quotient algebra $\mathcal{Q}/J$ has a unique CQG structure $(\mathcal{Q}/J, \Delta_0)$ such that the canonical quotient map $\pi_0: (\mathcal{Q}, \Delta) \to (\mathcal{Q}/J, \Delta_0)$ is a CQG morphism.
For any CQG morphism $\theta: (\mathcal{Q}_1, \Delta_1) \to (\mathcal{Q}_2, \Delta_2)$, the kernel $ker \theta$ is a Woronowicz $C^*$-ideal. If $J \subset ker \theta$ is a Woronowicz $C^*$-ideal, then there exists a unique CQG morphism $\bar{\theta}: \mathcal{Q}_1/J \to \mathcal{Q}_2$ such that $ \bar{\theta}\pi = \theta $, where $\pi$ denotes the quotient map from $\mathcal{Q}_1$ to $\mathcal{Q}_1/J$.
 
\begin{ex}  \label{Example_CQG} 
Now, let us review some important examples and constructions of CMQGs, which we will appear at later in this article.

\begin{enumerate}	
	\item Algebra of continuous function $C(G)$ on a compact Lie group $G \subseteq GL_{n}(\mathbb{C})$ \cite{Wor}, quantum permutation group $S_n^+ = (C(S_n^+),u)$ \cite{Wang}, and free unitary group $U_n^+ = (C(U_n^+),u)$ are some of the standard examples of CMQGs. The underlying $C^*$-algebra of the first one is commutative, whereas those of the latter two are non-commutative. It is worth noting that $C(S_n^+)/I_{com} \cong C(S_n)$ and $C(U_n^+)/I_{com} \cong C(U(n))$, where $I_{com}$ denotes the commutator ideal of respective algebras.\\
 	\item $C(H_{n}^{\infty +})$ (resp. $C(SH_{n}^{\infty +})$) is defined to be the universal $C^*$-algebra generated by $\{u_{ij}: i,j \in [n] \}$ such that
 \begin{itemize}
 	\item the matrices $u:=(u_{ij})_{n \times n}$ and $(u_{ij}^{*})_{n \times n}$ are unitaries. 
 	\item $u_{ij}$'s are normal partial isometries for all $i,j \in [n]$ (resp., $u_{ij}$'s are partial isometries for all $i,j \in [n]$).
 \end{itemize}
 	 Then $(C(H_{n}^{\infty +}), u)$ (resp. $(C(SH_{n}^{\infty +}), u)$) forms a CMQG, and it's denoted by $H_{n}^{+}$ (resp. $SH_n^{\infty +}$). In both cases, the quotient by the commutator ideal is isomorphic to $C(S^1 \wr S_n)$ (refer to \cite{Bichonwreath} for details). Also, both $H_n^{\infty +}$ and $SH_n^{\infty +}$ are quantum subgroups of $U_n^+$.\\
 	 \item For  $F \in \mathbb{GL}_{n}\mathbb{(C)}$, $ A_{u^{t}}(F)$ is the universal $C^{*}$-algebra generated by  $\{u_{ij}: i,j \in [n]\}$ subject to the conditions:
 	 \begin{itemize}
 	 	\item $u^{t}$ is unitary,
 	 	\item $uF^{-1}u^{*}F=F^{-1}u^{*}Fu=I_{n}$. 
 	 \end{itemize}
 	 Then $(A_{u^t}(F), u)$ forms a CMQG (see \cite{Van}, \cite{universalCQG} for details). Note that $(A_{u^t}(I_n), u) \cong U_n^+$.  \\

 	 \item Let $\{(A_i, q^{i})\}_{i \in [m]}$ be a finite family of CMQGs with multiplicative identities $\{1_{A_i}\}_{i \in [m]}$. Their {\bf free product} (denoted by $\mystar\limits_{i=1}^{m}  A_i$ or $A_1 * A_2 * \cdots * A_m$) is a CMQG $\left(\mystar\limits_{i=1}^{m} A_i, \oplus_{i=1}^{n} ~ q^{i} \right)$ \cite{Wangfree}, whose underlying $C^*$-algebra $\mystar\limits_{i=1}^{m} A_i$ is the unital universal $C^*$-algebra generated by $1$ and $\{q^{i}_{kl}\}$ such that $q^{i}_{kl}$'s fulfill the relations of $A_i$ for each $i \in [m]$ and $1=1_{A_1} = \cdots = 1_{A_m}$. Here the fundamental representation $\oplus_{i=1}^{n} ~ q^{i}$ denotes the block diagonal matrix $diag(q^{1}, q^{2}, ..., q^{m})$ whose diagonal blocks are $q^{1}, q^{2}, ..., q^{m}$ respectively.\\
 	 \item Let $(A,q)$ be a CMQG with fundamental representation $q$ of order $m \times m$. The {\bf free wreath product} (refer to \cite{Bichonwreath}, \cite{Free prod formulae} for details) of $(A,q)$ with $S_n^+ =(C(S_n^+),u)$ is a CMQG, simply denoted by $A \wr_{*} S_n^+$, whose underlying $C^*$-algebra is
 	 $\left(\mystar\limits_{k=1}^{n} A \right) * C(S_n^+) /<q_{ij}^{k}u_{kl}-u_{kl}q_{ij}^{k}>,$ where $q_{ij}^{k}$'s $(i,j \in [m])$ denote the generators of the $k$-th copy of $A$ in the free product $\left(\mystar\limits_{k=1}^{n} A, \oplus_{k \in [n]} ~q^{k} \right)$. The canonical fundamental representation $(q \wr_{*} u)$ is then given by the block matrix $(q^{k}u_{kl})_{k, l \in [n]}$, i.e. $(q \wr_{*} u)_{ik, jl} := q^{k}_{ij} u_{kl}$.\\
 	For example, it can be shown that $H_{n}^{\infty +} \cong C(S^1) \wr_{*} S_n^+$.
 \end{enumerate}	 
	
\end{ex}
Now, we prove the following elementary fact, which will be used later to prove Lemma \ref{lem_finitely many}.
\begin{lem} \label{A_u^t(F^Gamma)}
    Let $F=(f_{kl})_{n \times n}$ be a diagonal matrix $ diag(\underbrace{f_1,...,f_1}_{n_1-times}; \underbrace{f_2,..., f_2}_{n_2-times}; ...; \underbrace{f_p,..., f_p}_{n_p-times})$ such that $n=\sum_{i=1}^{p} n_i$ and $f_1>f_2>...>f_p$. Then $(A_{u^t}(F)/I_{com}, \bar{u})$ is a quantum subgroup of $\left(\otimes_{i=1}^{p}  C(U({n_i})), \oplus_{i=1}^{p} \bar{u}^{i} \right)$ and hence of $(C(U(n)),u)$, where $I_{com}$ denotes the commutator ideal of $A_{u^t}(F)$ and $\bar{u}$ denotes the image of $u$ under the canonical quotient map $A_{u^t}(F) \to A_{u^t}(F)/I_{com}$.
\end{lem}
\begin{proof}
     Since $\bar{u}^*F\bar{u}=F$ and ${{\bar{u}}^{t}} {\bar{u}^t}{}^{*}=1$, we have 
     \begin{equation} \label{Au^tF_eq1}
         \sum\limits_{i=1}^{n} \bar{u}_{ki}^* f_{kk} \bar{u}_{ki} = f_{ii}.1 = f_{ii} (\sum\limits_{k=1}^{n}  \bar{u}_{ki} \bar{u}_{ki}^{*})=f_{ii} (\sum\limits_{k=1}^{n}  \bar{u}_{ki}^{*} \bar{u}_{ki}).
     \end{equation} 
     For $i=1,2,...,n_1$, Equation \eqref{Au^tF_eq1} implies 
     \begin{align*}
         & \sum_{k=1}^{n} (f_{ii} - f_{kk}) \bar{u}_{ki}^* \bar{u}_{ki} = 0 \\
\implies & \sum_{k=n_1 +1}^{n} (f_{ii} - f_{kk}) \bar{u}_{ki}^* \bar{u}_{ki}=0. 
     \end{align*}
Since $f_1> f_j$ for all $j=2, 3, ..., p$, the above equation implies $\bar{u}_{ki}=0$ for all $k \in [n]-[n_1]$ and $i=1, 2, ..., n_1$. Since $F$ is diagonal, applying $\kappa$ to the last equation, we further have $\bar{u}_{ik}=0$ for the same indices $i,k$.  \\
Similarly, for $i=n_1+1, ..., n_1 + n_2$, the above fact and Equation \eqref{Au^tF_eq1} give 
    \begin{align*}
          & \sum_{k=n_1 +1}^{n} (f_{ii} - f_{kk}) \bar{u}_{ki}^* \bar{u}_{ki} = \sum_{k=n_1+n_2+1}^{n} (f_{ii} - f_{kk}) \bar{u}_{ki}^* \bar{u}_{ki} = 0
     \end{align*}
Hence, $\bar{u}_{ki}=0=\bar{u}_{ik}$ for all $k \in [n]-[n_1+n_2]$ and $i=n_1+1, n_1+2, ..., n_2$. Continuing this argument, $\bar{u}$ reduces to the block diagonal matrix $\bar{u}^{i}$ of each block of order $n_i \times n_i$ for each $i \in [p]$. Now, relation $\bar{u}F^{-1}\bar{u}^{*}F=F^{-1}\bar{u}^{*}F\bar{u}=I_{n}$ actually implies $\bar{u}$ is unitary. Thus, $(A_{u^t}(F^{\Gamma})/I_{com},\bar{u})$ is a quantum subgroup of $(C(U(n)),u)$, and more precisely, it is a quantum subgroup of $\left(\otimes_{i=1}^{p}  C(U({n_i})), \oplus_{i=1}^{p} \bar{u}^{i} \right)$.
\end{proof}
\subsection{Quantum automorphism group of a graph}
For the purposes of this article, we recall Banica's notion of the quantum automorphism group of a directed graph without multiple edges. \\ 
Let $\Gamma=\{V(\Gamma),E(\Gamma), s, r \}$ be a finite directed graph without multiple edges on $n$ vertices with  $V(\Gamma)=\{1,2,...,n\}$. The {\bf quantum automorphism group} of $\Gamma$ (in the sense of Banica), denoted by $QAut_{Ban}(\Gamma)$, is the CMQG whose underlying $C^*$-algebra is the universal $C^{*}$-algebra with generators $u_{ij}$ $(i, j \in [n])$ satisfying the relations
(i) $u_{ij}=u_{ij}^{*}=u_{ij}^{2}$, 
(ii)	$\sum_{l=1}^{n}u_{il}=1=\sum_{l=1}^{n}u_{li}$, 
(iii) $u A_{\Gamma}=A_{\Gamma} u$,
and $u=(u_{ij})_{n \times n}$ is the (canonical) fundamental representation.
If no ambiguity arises, we shall use the notation $QAut_{Ban}(\Gamma)$ for both the CMQG and its underlying $C^*$-algebra.
Relation (iii) is equivalent to the following relations [Lemma 3.1.1, \cite{Fulton}]:
\begin{align}
   & u_{ij}u_{kl}=0=u_{kl}u_{ij} \text{ if } i \rightarrow k \text{ and } (j,l) \notin E(\Gamma), \label{adj_nonadj_prod=0_eq1}\\
& u_{ij}u_{kl}=0=u_{kl}u_{ij} \text{ if } (i,k) \notin E(\Gamma) \text{ and } j \rightarrow l. \label{adj_nonadj_prod=0_eq2} 
\end{align}
 
\noindent [Theorem 2.1.6, \cite{SSthesis}] ensures that the quantum automorphism group of a graph always contains its classical automorphism group in the following sense:\\	
If $A:=QAut_{Ban}(\Gamma)/<u_{ij}u_{kl}-u_{kl}u_{ij}>$ is the quotient of $QAut_{Ban}(\Gamma)$ by the commutator ideal, then $A \cong C(Aut(\Gamma))$, where $Aut(\Gamma)$ denotes the classical automorphism group of $\Gamma$.\\
Thus, in the view of the above, there exist surjective *-homomorphisms
$$C(S_n^+) \to QAut_{Ban}(\Gamma) \to C(Aut(\Gamma))$$ sending the canonical generators to the generators. We say the $\Gamma$ is asymmetric (resp. quantum asymmetric) if $Aut(\Gamma) \cong \{e\}$ (resp. $QAut_{Ban}(\Gamma) \cong \mathbb{C}$).\\[0.2cm]

\noindent Next, we recall a few results that will be used later in this article.\\[0.1cm]
{\bf [Lemma 3.2.3, \cite{Fulton}]:} Let $u$ be the canonical fundamental representation of ${QAut}_{Ban}(\Gamma)$. If $indeg(i) \neq indeg(j)$ or $outdeg(i) \neq outdeg(j)$ for some $i,j \in V(\Gamma)$, then $u_{ij} = 0$.	\\

\noindent Two graphs $\Gamma_1$ and $\Gamma_2$ are said to be {\bf quantum isomorphic} if there exists a $C^*$-algebra $\mathcal{A}$ with generators $u_{ij}$'s (where $i \in V(\Gamma_1)$ and $j \in V(\Gamma_2)$) satisfying the relations: (1) $u_{ij}=u_{ij}^{*}=u_{ij}^{2}$, (2) $\sum_{i \in V(\Gamma_1)}u_{ij}=1=\sum_{j \in V(\Gamma_2)}u_{ij}$, (3) $u A_{\Gamma_1}=A_{\Gamma_2} u$. Clearly, in this case, $|V(\Gamma_1)|=|V(\Gamma_2)|$.\\

We recall the behavior of the quantum symmetry of a disjoint union of connected graphs in the setting of simple undirected graphs.\\[0.2cm]
{\bf [Lemma 6.3, \cite{Qauttrees}]:} Let $\Gamma$ be a simple undirected graph with connected components $\{\Gamma_{i}\}_{i=1}^{m}$. If $u_{pq} \neq 0$ for some $p \in V(\Gamma_i)$ and $q \in V(\Gamma_j)$ for $i \neq j$, then $\Gamma_{i}$ and $\Gamma_{j}$ are quantum isomorphic.\\[0.2cm]
{\bf [Lemma 6.4, \cite{Qauttrees}]:} Let $\{\Gamma_{i}\}_{i=1}^{m}$ be some simple undirected graphs such that each connected component of $\Gamma_{i}$ is quantum non-isomorphic to a connected component of $\Gamma_{j}$ for all $i \neq j$. Then, $ QAut_{Ban}(\sqcup_{i=1}^{m} \Gamma_{i}) \cong \mystar\limits_{i=1}^{m} QAut_{Ban}(\Gamma_{i}) $.\\

[Theorem 4.5, \cite{Lupini}] \& [Corollary 7.1.4, \cite{SSthesis}] also establishes the similar result for $m=2$. Later, in Proposition \ref{graph_free prod}, we ensure that similar results hold in the setting of directed graphs.

\subsection{Quantum symmetry of graph $C^*$-algebra}
We first review the notion of quantum symmetry of a graph $C^*$-algebra introduced by Joardar and Mandal \cite{Mandal}. We consider the finite directed graphs $\Gamma=\{V(\Gamma),E(\Gamma),s,r\}$ without isolated vertices. Since $\Gamma$ has no isolated vertex, it is enough to define an action on the partial isometries $\{S_e\}_{e \in E(\Gamma)}$ corresponding to its edges.\\ 
A faithful action $\alpha$ of a CQG $(Q, \Delta)$ on a graph $C^{*}$-algebra $C^{*}(\Gamma)$ is said to be {\bf linear} if $ \alpha(S_{e})=\sum\limits_{f \in E(\Gamma)} S_{f} \otimes q_{fe}$, where $q_{ef}\in Q$ for each $e,f \in E(\Gamma)$.\\
Consider the sets:
\begin{itemize}
	\item[(1)] $\mathcal{I}:=\{u\in V(\Gamma) : s^{-1}(u) = \emptyset\}$, 
	\item[(2)] $ E':= \{(e,f) \in E(\Gamma) \times E(\Gamma) : S_{e}S_{f}^{*} \neq 0 \}=\{(e,f) \in E(\Gamma) \times E(\Gamma) : r(e)=r(f) \}$.
\end{itemize}

\noindent Now, Lemma 3.2 of \cite{Mandal} ensures that $\{ p_{u}, S_{e}S_{f}^{*} : u \in \mathcal{I}, (e,f) \in E' \}$ is a linearly independent set. Thus, consider the subspace $\mathcal{V}_{2,+}= span\{ p_{u}, S_{e}S_{f}^{*} : u \in \mathcal{I}, (e,f) \in E' \}$ and define a linear functional $\tau: \mathcal{V}_{2,+} \to \mathbb{C}$ by $\tau(S_{e}S_{f}^{*})=\delta_{ef}$, $\tau(p_{u})=1 $ for all $(e,f) \in E' $ and $ u \in \mathcal{I}$
(see subsection 3.1 of \cite{Mandal}).\\[0.1cm]
Since $\alpha(\mathcal{V}_{2,+}) \subseteq \mathcal{V}_{2,+} \otimes Q$ (Lemma 3.6 of \cite{Mandal}), the equation $(\tau \otimes id)\alpha(.)=\tau(.).1 $ on $\mathcal{V}_{2,+}$ makes sense, and we say that $\alpha$ preserves the functional $\tau$.\\
 Now, define a {\bf category $\mathfrak{C}_{\tau}^{Lin} $} [Definition 3.7, \cite{Mandal}] whose objects are $(Q,\Delta),\alpha)$ acting on $C^*(\Gamma)$ via a linear faithful such that $\alpha$ preserves $\tau$ on $\mathcal{V}_{2,+}$. A morphism from $((Q_1,\Delta_{1}),\alpha_{1})$ to $((Q_2,\Delta_{2}),\alpha_{2})$ is given by a CQG morphism $\phi:Q_1 \to Q_2 $ such that $ (id_{C^{*}(\Gamma)} \otimes \phi)\alpha_{1}=\alpha_{2} $.\\

\noindent Define the matrix $F^{\Gamma}$ of order $|E(\Gamma)| \times |E(\Gamma)|$ given by $(F^{\Gamma})_{ef}= \tau(S_{e}^{*}S_{f})$. Equivalently, by [Lemma 4.3, \cite{Mandal}],  $F^{\Gamma}$ is an invertible diagonal matrix with $$(F^{\Gamma})_{ee}= \begin{cases}
	|s^{-1}(r(e))| & if ~ s^{-1}(r(e)) \neq \emptyset\\
	1 & if ~ s^{-1}(r(e)) = \emptyset
	\end{cases}.$$
Thus, $A_{u^{t}}(F^{\Gamma})$ is a CMQG. 
It has been shown in [Proposition 3.8, \cite{Mandal}] that there is a surjective CQG morphism from $A_{u^{t}}(F^{\Gamma})$ to any object in category $\mathfrak{C}_{\tau}^{Lin}$, and consequently standard arguments show that the category $\mathfrak{C}_{\tau}^{Lin} $ admits a {\bf universal object} [Theorem 3.9, \cite{Mandal}], denoted by $((Q_{\tau}^{Lin}(\Gamma), \Delta), \alpha)$, in the following sense:\\
Given any object $((Q_0, \Delta_0),\alpha_0)$ of the category $\mathfrak{C}_{\tau}^{Lin} $, there is a surjective CQG morphism $\phi: (Q_{\tau}^{Lin}(\Gamma), \Delta) \to (Q_0, \Delta_0) $ such that $(id_{C^*(\Gamma)} \otimes  \phi)\alpha=\alpha_{0} $.\\
Since the finite set $\{S_e\}_{e \in E(\Gamma)}$ generates $C^*(\Gamma)$, the underlying CQG $(Q_{\tau}^{Lin}(\Gamma), \Delta)$ of the universal object is, in fact, a CMQG with the fundamental representation $q$ of order $|E(\Gamma)| \times |E(\Gamma)|$ determined by the action $\alpha(S_e)=\sum\limits_{f \in E(\Gamma)} S_f \otimes q_{fe}$.\\
We refer to the CQG $(Q_{\tau}^{Lin}(\Gamma), \Delta)$, or equivalently the CMQG $(Q_{\tau}^{Lin}(\Gamma), q)$ (or simply $Q_{\tau}^{Lin}(\Gamma)$ if no confusion arises), as the {\bf quantum symmetry of $C^*(\Gamma)$}. The matrix $q$ is called the canonical fundamental representation of $Q_{\tau}^{Lin}(\Gamma)$.
Note in particular that $(Q_{\tau}^{Lin}(\Gamma),q)$ is a quantum subgroup of $(A_{u^t}(F^{\Gamma}), u)$ in the sense that the map $\phi: A_{u^t}(F^{\Gamma}) \to  Q_{\tau}^{Lin}(\Gamma)$ sends  $u_{ij} \mapsto q_{ij}$.\\
Further, we have the following:\\[0.1cm]
$\bullet$ {\bf [Proposition 3.11, \cite{Mandal}]:} There exists a surjective CQG morphism $Q_{\tau}^{Lin}(\Gamma) \to \mystar\limits_{|E(\Gamma)|} C(S^1)$ maps $q_{ee} \mapsto u_{e}$ and $q_{ef} \mapsto 0$ for $e \neq f$, where $\{u_{e}\}_{e \in E(\Gamma)}$ are unitaries generating $\left(\mystar\limits_{|E(\Gamma)|} C(S^1),  \Delta\right)$.\\
Consequently, the $C^{*}$-algebra of $Q_{\tau}^{Lin}(\Gamma)$ is noncomutative for a graph $\Gamma $ with $|E(\Gamma)| \geq 2 $ [Corollary 3.12, \cite{Mandal}].\\[0.1cm]
$\bullet$ {\bf [Theorem 3.10, \cite{Mandal}]:} $((QAut_{Ban}(\Gamma), \Delta), \alpha)$ is an object of the category $\mathfrak{C}_{\tau}^{Lin}$, where the action $\alpha: C^*(\Gamma) \to C^*(\Gamma) \otimes QAut_{Ban}(\Gamma)$ is given by 
	\begin{align*}
	\alpha(p_i) &= \sum\limits_{k \in V(\Gamma)} p_k \otimes u_{ik} \\
	\alpha(S_e) &= \sum\limits_{f \in E(\Gamma)} S_{f} \otimes u_{s(e)s(f)} u_{r(e)r(f)}.
	\end{align*}
    \vspace{0.2cm}


Moreover, we recall several families of graphs $\Gamma$ for which the quantum symmetry of the associated graph $C^*$-algebra $C^*(\Gamma)$ has already been identified. \\[0.1cm]
\noindent $\bullet$ {\bf [Theorem IV.6, \cite{rigidity}]:} Consider the family $\mathcal{G}$ (as described in [Section III, \cite{rigidity}]) of graphs $\Gamma=\{V(\Gamma),E(\Gamma), s, r\}$ that satisfy the following properties:
\textbf{(R1)} There does not exist any cycle of length $\geq$ 2;
\textbf{(R2)} There exists a path of length $(|V(\Gamma)|-1)$ that consists of all the vertices;
\textbf{(R3)} Given any two vertices (which may or may not be distinct), there exists at most one edge joining them.\\
Then for any $\Gamma \in \mathcal{G} $, $ Q_{\tau}^{Lin}(\Gamma) \cong \left(\mystar\limits_{|E(\Gamma)|} C(S^1),  \Delta\right)$.\\[0.1cm]
In Examples 1--5 of Section III of \cite{rigidity}, the authors discuss several well-known classes of graphs that belong to the family $\mathcal{G}$.\\[0.2cm]
$\bullet$ {\bf [Proposition 3.3, \cite{max_permutational_symm}]:} Let $\mathcal{S}$ be a family of graphs $\Gamma$ without loops such that $r(e)=v$ for all $e \in E(\Gamma)$ for some fixed vertex $v \in V(\Gamma)$. If $\Gamma \in \mathcal{S}$, then $Q_{\tau}^{Lin}(\Gamma) \cong U_{|E(\Gamma)|}^+$.\\[0.1cm]
$\bullet$ {\bf [Proposition 4.2, \cite{Mandalkms}]:} $Q_{\tau}^{Lin}(\sqcup_{i=1}^{n} L_1) \cong H_{n}^{\infty+}$.\\[0.1cm]
$\bullet$ {\bf [Theorem 4.6, \cite{Mandal}]:} $Q_{\tau}^{Lin}(C_2) \cong H_{2}^{\infty+}$.\\[0.1cm]
$\bullet$ {\bf [Proposition 3.4, \cite{max_permutational_symm}]:} Let $\mathcal{I}_1$ be a family of graphs $\Gamma$ without loops such that for each $e \in E(\Gamma) $, $r(e)$ is a sink of indegree 1. If $\Gamma \in \mathcal{I}_1$, then $Q_{\tau}^{Lin}(\Gamma) \cong SH_{|E(\Gamma)|}^{\infty+}$.

\section{Main Section}
In this section, we highlight several characteristic differences between the nature of the quantum symmetry group of a directed graph and that of its associated graph $C^*$-algebra. All the graphs considered here are {\bf finite directed graphs without isolated vertices and without multiple edges (that is, whose adjacency matrices are $\{0,1\}$-matrices such that, for each index $i$, the $i^{th}$ row and $i^{th}$ column are not both zero vectors)}, so that one can systematically compare the quantum symmetry groups of graphs in the sense of Banica as well as the quantum symmetry groups of their graph $C^*$-algebras in the sense of Joardar and Mandal within a common framework. We divide this section into a few subsections, each of which reflects a comparison of the nature of quantum symmetry groups in the two frameworks.\\

Let $\mathcal{G}_{0,1}$ denote the family of finite directed graphs $\Gamma$ without isolated vertices and without multiple edges (but possibly with loops). We define $\mathcal{Q}_{Ban}:= \{QAut_{Ban}(\Gamma): \Gamma \in \mathcal{G}_{0,1} \}$,
$\mathcal{Q}_{JM}:= \{Q_{\tau}^{Lin}(\Gamma): \Gamma \in \mathcal{G}_{0,1} \}$.

\subsection{No common quantum symmetry} \label{sec_no common q. symm}
It is known that, for all such graphs $\Gamma$ mentioned above, $QAut_{Ban}(\Gamma)$ is a quantum subgroup of $S_{|V(\Gamma)|}^+$ and $Q_{\tau}^{Lin}(\Gamma)$ is a quantum subgroup of $A_{u^t}(F^{\Gamma})$ for some invertible scalar matrix $F^{\Gamma}$ of order $|E(\Gamma)|$. Here, we ensure that there are no common CMQGs that can appear as the quantum symmetry group for both graph and graph $C^*$-algebraic settings.

\begin{prop} \label{no common q. symm.}
	Let $Q$ be a CMQG such that $Q \cong QAut_{Ban}(\Gamma)$ for some graph $\Gamma$ without multiple edges. Then there does not exist any graph $\Gamma'$ without isolated vertices for which $Q_{\tau}^{Lin}(\Gamma') \cong Q$. In particular, $\mathcal{Q}_{Ban} \cap \mathcal{Q}_{JM} = \emptyset$.
\end{prop}
\begin{proof}
	If possible, there exists a graph $\Gamma'$ (without isolated vertices) for which $Q_{\tau}^{Lin}(\Gamma') \cong Q \cong QAut_{Ban}(\Gamma)$. By [Proposition 3.11, \cite{Mandal}], there always exists a surjective $C^*$-homomorphism from $Q_{\tau}^{Lin}(\Gamma')$ to  $\mystar\limits_{|E(\Gamma')|} C(S^1)$. Thus, there exists a surjective $C^*$-homomorphism $\phi: QAut_{Ban}(\Gamma)\to C(S^1)$ and hence a surjective $C^*$-homomorphism from $C(S_{|V(\Gamma)|}^+)$ to $C(S^1)$. Let $I_{com}$ denote the commutator ideal of $C(S_{|V(\Gamma)|}^+)$. Since $C(S_{|V(\Gamma)|}^+)/I_{com} \cong C(S_{|V(\Gamma)|})$, we again obtain a surjective $C^*$-homomorphism from $C(S_{|V(\Gamma)|})$ to $C(S^1)$, a contradiction.  
\end{proof}
\vspace{0.4cm}
\subsection{Same quantum symmetry appears for infinitely many graphs but for finitely many graph $C^*$-algebras} \label{sec_q. symm finite grp C*-alg}
Given a CMQG $Q$ in $\mathcal{Q}_{Ban}$, there are infinitely many non-isomorphic graphs in $\mathcal{G}_{0,1}$ that attain the same quantum symmetry group $Q$ (Theorem \ref{graph_infinitely many}); whereas in the graph $C^*$-algebraic scenario, for a given $Q$ in $\mathcal{Q}_{JM}$, there are only finitely many graphs in $\mathcal{G}_{0,1}$ whose associated graph $C^*$-algebra attains the quantum symmetry $Q$ (Theorem \ref{graph C*alg_finitely many}).\\

\noindent $\blacksquare$ For the graph scenario, we already established a similar result in [Corollary 5.4 \& Remark 5.5,  \cite{some results bichon_HKM}] for undirected graphs. We encode the same idea to do it for a directed graph. The idea of the proof is exactly the same as in [Lemma 5.1, 5.2, \cite{some results bichon_HKM}]. We only modify the construction analogously for the directed case and refer to the necessary results from Chapter 3 of \cite{Fulton} used here. We state and briefly sketch the proof for the directed scenario in the following theorem.
\begin{thm}\label{graph_infinitely many}
	Given a directed graph $\Gamma$ (without multiple edges) on $n (>1)$ vertices, there always exists an infinite family of mutually non-isomorphic weakly connected directed graphs  $\{\Gamma_i\}_{i \in \mathbb{N}}$ such that $QAut_{Ban}(\Gamma_i) \cong QAut_{Ban}(\Gamma)$ for all $i \in \mathbb{N}$.
\end{thm}
\begin{proof}
We divide the proof into three steps.\\
	{\bf STEP-1:} We claim that if $\Gamma$ is not weakly connected, then there always exists a weakly connected graph $\Gamma'$ such that $QAut_{Ban}(\Gamma) \cong QAut_{Ban}(\Gamma')$.\\
	 To that end, construct $\Gamma'$ by adjoining a new vertex $v \notin V(\Gamma)$ as follows:
	$$V(\Gamma'):= V(\Gamma) \sqcup \{v\} \text{ \ \ \ and \ \ \ } E(\Gamma'):= E(\Gamma) \sqcup \{(v,w), (u,v): w,u \in V(\Gamma) \}.$$ 
	Then clearly $A_{\Gamma'}= \begin{bmatrix}
		0 & J_{1,n} \\
		J_{n,1} & A_{\Gamma}
	\end{bmatrix} $.
	Let $u$ and $u'$ be the canonical fundamental representations of $QAut_{Ban}(\Gamma)$ and $QAut_{Ban}(\Gamma')$ respectively. Since $\Gamma$ is disconnected, $outdeg_{\Gamma'}(v)=indeg_{\Gamma'}(v)=|V(\Gamma)| > indeg_{\Gamma'}(a)$ for all $a \in V(\Gamma)$. It follows from [Lemma 3.2.2, \cite{Fulton}] that  $u'_{vv}=1$. Hence $u'$ reduces to the block form $\begin{bmatrix}
	1 & 0_{1,n}\\
	0_{n,1} & u''
	\end{bmatrix}$, where $u''=(u'_{ab})_{a,b \in V(\Gamma)}$ is also a magic unitary matrix. It is easy to observe that $u' A_{\Gamma'}=A_{\Gamma'}u'$ is equivalent to $u'' A_{\Gamma}=A_{\Gamma}u''$. Thus, using the universal properties of the respective algebras, we obtain our claim. \\
	
	{\bf STEP-2:} Let $\Gamma'$ be a weakly connected graph with $k$ vertices, say $\{v_i\}_{i=1}^k$. Construct a weakly graph $\Gamma_1$ from $\Gamma'$ by adjoining $k$ new vertices (say $\{v_{k+i}\}_{i=1}^{k}$) as follows: 
	$$V(\Gamma_1):= V(\Gamma') \sqcup \{v_{k+i}: i \in [k]\} \text{ \ \ \ and \ \ \ } E(\Gamma_1):= E(\Gamma') \sqcup \{(v_i,v_{k+i}), (v_{k+i}, v_i ): i \in [k] \}.$$
	First, note that $A_{\Gamma_1} = \begin{bmatrix}
	A_{\Gamma} & I_k  \\
	I_k & 0_k
	\end{bmatrix}.$
	Let $u_1$ denote the canonical fundamental representation of $QAut_{Ban}(\Gamma_1)$.
	Observe that at least one of $indeg_{\Gamma_1}(v_i)$ or $outdeg_{\Gamma_1}(v_i)$ is strictly greater than 1 (since $\Gamma'$ is weakly connected), however, $indeg_{\Gamma_1}(v_{k+i})=outdeg_{\Gamma_1}(v_{k+i}) =1$ for all $i \in [n]$. Thus, using [Lemma 3.2.2, \cite{Fulton}], the matrix $u_1$ can be decomposed into block form $\begin{bmatrix}
	u_1' & 0_{k}\\
	0_{k} & u_1''
	\end{bmatrix}$, where $u_1'$ and $u_1''$ are two magic unitary matrices of order $k$.
	A direct computation shows that the relation $u_1 A_{\Gamma_1}=A_{\Gamma_1}u_1$ is equivalent to the relations $u_1' A_{\Gamma}=A_{\Gamma} u_1'$ and $u_1' = u_1''$. Now, the standard argument shows that $QAut_{Ban}(\Gamma') \cong QAut_{Ban}(\Gamma_1)$ (for more details refer to [Lemma 5.2, \cite{Fulton}]).\\
	
	{\bf STEP-3:} Given a graph $\Gamma$, using STEP-1 and STEP-2, we can always produce a weakly connected graph $\Gamma_1$ having distinct order such that $QAut_{Ban}(\Gamma) \cong QAut_{Ban}(\Gamma_1)$. Now, repeated application of STEP-2 yields an infinite family of graphs $\{\Gamma_i\}_{i \in \mathbb{N}}$ having mutually distinct orders (hence, mutually non-isomorphic) such that $QAut_{Ban}(\Gamma_i) \cong QAut_{Ban}(\Gamma)$ for all $i \in \mathbb{N}$.
\end{proof}
\vspace{0.4cm}
\noindent $\blacktriangleright$ On the other hand, in the graph $C^*$-algebraic setting, the analogue of Theorem \ref{graph_infinitely many} is not true. For instance, $C(S^1) \in \mathcal{Q}_{JM}$ can appear as the quantum symmetry of $C^*{(P_1)}$ and $C^*(L_1)$ only, because $Q_{\tau}^{Lin}(\Gamma)$ is non-commutative if $|E(\Gamma)|>1$ (by Corollary 3.12 of \cite{Mandal}). Thus, here we have only finitely many choices of graphs. We will show in Theorem \ref{graph C*alg_finitely many} that this is actually true in general as a consequence of the following lemma.

\begin{lem} \label{lem_finitely many}
If $\Gamma$ is a directed graph without isolated vertices, then $\mystar\limits_{m} C(S^1)$ cannot be a quantum subgroup of $Q_{\tau}^{Lin}(\Gamma)$ for $m > |E(\Gamma)|$.
\end{lem}
\begin{proof}
	Let $n:=|E(\Gamma)|$. Let also assume that there exists a surjective CQG morphism $\phi: Q_{\tau}^{Lin}(\Gamma) \to \mystar\limits_{m} C(S^1)$. Then, using [Proposition 3.8, \cite{Mandal}], we have the following sequence of surjective CQG morphisms:
    $$A_{u^t}(F^{\Gamma}) \to Q_{\tau}^{Lin}(\Gamma) \to \mystar\limits_{m} C(S^1) \to \otimes_{m} C(S^1). $$ Taking this all together, there exists a surjective CQG morphism from $\phi_1: A_{u^t}(F^{\Gamma}) \to \otimes_{m} C(S^1)$. Let $I_{com}:=<u_{ij}u_{kl}-u_{kl}u_{ij}>$ denote the commutator ideal of $A_{u^t}(F^{\Gamma})$. Since $I_{com} \subset ker\phi_{1}$, by [Theorem 2.11 (3),\cite{Wangfree}], there is a surjective CQG morphism from $A_{u^t}(F^{\Gamma})/I_{com}$ to $\otimes_{m} C(S^1)$. Consequently, by Lemma \ref{A_u^t(F^Gamma)}, there is a surjective CQG morphism from $C(U(n))$ to $ C((S^1)^m)$.  Hence, there exists a continuous injective group homomorphism from $(S^1)^m$ to $U(n)$ for $m>n$, which is not possible. 
\end{proof}

\begin{thm} \label{graph C*alg_finitely many}
Given a graph $\Gamma$ without isolated vertices, there exist only finitely many graphs $\{\Gamma_i\}_{i \in \Lambda}$ without isolated vertices (where $\Lambda$ is a finite index set) such that $Q_{\tau}^{Lin}(\Gamma_i) \cong Q_{\tau}^{Lin}(\Gamma)$.	
\end{thm}
\begin{proof}
	For a given graph $\Gamma$, if there exists a graph $\Gamma'$ such that  $Q_{\tau}^{Lin}(\Gamma') \cong Q_{\tau}^{Lin}(\Gamma)$ then $\mystar\limits_{|E(\Gamma)|} C(S^1)$ is a quantum subgroup of $Q_{\tau}^{Lin}(\Gamma')$. By Lemma \ref{lem_finitely many}, $|E(\Gamma')| \leq |E(\Gamma)| $. Since we consider graphs $\Gamma'$ without isolated vertices, there are at most finitely many choices.
\end{proof}
An immediate consequence of Lemma \ref{lem_finitely many} is the following:
\begin{cor} \label{edge equal}
If $Q_{\tau}^{Lin}(\Gamma_1) \cong Q_{\tau}^{Lin}(\Gamma_2)$, then $|E(\Gamma_1)|=|E(\Gamma_2)|$.
\end{cor}
\begin{proof}
Since there is a surjective CQG morphism from $Q_{\tau}^{Lin}(\Gamma_2) \to \mystar\limits_{|E(\Gamma_2)|} C(S^1)$ (by [Proposition 3.11, \cite{Mandal}]), and given that $Q_{\tau}^{Lin}(\Gamma_1) \cong Q_{\tau}^{Lin}(\Gamma_2)$, there exists a surjective CQG morphism from $Q_{\tau}^{Lin}(\Gamma_1)$ to $\mystar\limits_{|E(\Gamma_2)|} C(S^1)$ (i.e., $\mystar\limits_{|E(\Gamma_2)|} C(S^1)$ is a quantum subgroup of $Q_{\tau}^{Lin}(\Gamma_1)$). By Lemma \ref{lem_finitely many}, $|E(\Gamma_2)| \leq |E(\Gamma_1)|$. A similar argument shows that $|E(\Gamma_1)| \leq |E(\Gamma_2)|$, and hence the result follows.
\end{proof}
The above property shows that the behaviour of the quantum symmetry group of a graph $C^*$-algebra is more rigid than that of the quantum symmetry group of a graph, which also leads to a difference:

\begin{rem}
By construction in Theorem \ref{graph_infinitely many}, we have the following: Given $Q \in \mathcal{Q}_{Ban}$, we can construct a graph with an arbitrary large cardinality of vertex set and edge set having quantum symmetry group $Q$. \\
In contrast, if $Q \in \mathcal{Q}_{JM}$ with $Q \cong Q_{\tau}^{Lin}(\Gamma_0)$, then we cannot find a graph $\Gamma$ with $|E(\Gamma)| > |E(\Gamma_0)|$ such that $Q_{\tau}^{Lin}(\Gamma) \cong Q$.  
\end{rem}

\subsubsection{Uniqueness of quantum symmetry}
(1) It may happen that there exists a quantum symmetry group $Q \in \mathcal{Q}_{JM}$ for which there is a unique graph $\Gamma_0 \in \mathcal{G}_{0,1}$ such that $Q_{\tau}^{Lin}(\Gamma_0) \cong Q$. For example, the free unitary group $U_2^+$ appears as $Q_{\tau}^{Lin}(Si_2)$ for the only choice of graph $Si_2 \in \mathcal{G}_{0,1}$ (see FIGURE \ref{fig Si_2}). 
\begin{prop} \label{prop U_2^+}
	If $Q_{\tau}^{Lin}(\Gamma) \cong U_2^+$ for some $\Gamma \in \mathcal{G}_{0,1}$, then $\Gamma \cong Si_2$.
\end{prop}
\begin{proof}
	Using Corollary \ref{edge equal}, first observe that if $Q_{\tau}^{Lin}(\Gamma) \cong U_2^+$, then $|E(\Gamma)|=2$. Thus, we can identify all possible graphs $\Gamma \in \mathcal{G}_{0,1}$ with exactly 2 edges and their quantum symmetry groups $Q_{\tau}^{Lin}(\Gamma)$, as shown in the table bellow:\\
	\begin{center}
		\begin{tabular}{c|c|c}
			\hline
			 $\Gamma \in \mathcal{G}_{0,1}$ & $Q_{\tau}^{Lin}(\Gamma)$ & using the result \\
			\hline
     
       \begin{minipage}{0.3\linewidth}
       	\centering
       		 	\begin{tikzpicture}   
       		\draw[fill=black] (0,0) circle (2pt);
       		\draw[fill=black] (2,0) circle (2pt);
       		\draw[fill=black] (4,0) circle (2pt);

       		\draw[black] (0,0)-- node{\rmidarrow}  (2,0);
       		\draw[black] (2,0)-- node{\rmidarrow} (4,0);
       		\end{tikzpicture} \captionof{figure}{$P_2$}
       \end{minipage}
        & $C(S^1)*C(S^1)$ & [Theorem IV.6, \cite{rigidity}]  \\
		 \hline
		\begin{minipage}{0.3\linewidth}
			\centering
				\begin{tikzpicture}
			\draw[fill=black] (0,0) circle (2pt);
			\draw[fill=black] (2,0) circle (2pt);
			
			\begin{scope}[decoration={markings,
				mark=at position 0.25 with {\arrow[black,thick]{<}},
			}]
			\draw[black,postaction={decorate}, thick](0,0.5) circle (0.5);
			\end{scope}
			
			\begin{scope}[decoration={markings,
				mark=at position 0.50 with {\arrow[black,thick]{>}},
			}]
			\draw[black,postaction={decorate}, thick ](0,0)--(2,0);
			\end{scope}		
			\end{tikzpicture}   \captionof{figure}{$T$}
		\end{minipage}		
	  & $C(S^1)*C(S^1)$ &  [Theorem IV.6, \cite{rigidity}]  \\
			\hline
			\begin{minipage}{0.3\linewidth}
				\centering
				\begin{tikzpicture}
				\draw[fill=black] (0,0) circle (2pt);
				\draw[fill=black] (2,0) circle (2pt);
				
				\begin{scope}[decoration={markings,
					mark=at position 0.25 with {\arrow[black,thick]{<}},
				}]
				\draw[black,postaction={decorate}, thick](0,0.5) circle (0.5);
				\end{scope}
				
				\begin{scope}[decoration={markings,
					mark=at position 0.50 with {\arrow[black,thick]{<}},
				}]
				\draw[black,postaction={decorate}, thick ](0,0)--(2,0);
				\end{scope}		
				\end{tikzpicture} \captionof{figure}{$T'$}
			\end{minipage}
   	 & $C(S^1)*C(S^1)$ & [Theorem IV.6, \cite{rigidity}]  \\
			\hline
			\begin{minipage}{0.3\linewidth}
				\centering
				\begin{tikzpicture}
				\draw[fill=black] (0,0) circle (2pt);
				\draw[fill=black] (2,0) circle (2pt);
				\draw[fill=black] (3,0) circle (2pt);
				
				\begin{scope}[decoration={markings,
					mark=at position 0.25 with {\arrow[black,thick]{<}},
				}]
				\draw[black,postaction={decorate}, thick](3,0.5) circle (0.5);
				\end{scope}
				
				\begin{scope}[decoration={markings,
					mark=at position 0.50 with {\arrow[black,thick]{>}},
				}]
				\draw[black,postaction={decorate}, thick ](0,0)--(2,0);
				\end{scope}	
				\end{tikzpicture}  \captionof{figure}{$P_1 \sqcup L_1$}
			\end{minipage}
			 & $C(S^1)*C(S^1)$ & [Lemma IV.1, \cite{rigidity}]  \\
			\hline		 
		\begin{minipage}{0.3\linewidth}
			\centering
			\begin{tikzpicture}
			\draw[fill=black] (0,0) circle (2pt);
			\draw[fill=black] (2,0) circle (2pt);

			\draw[black,thick](0,0) edge[bend left= 60] node{\rmidarrow} (2,0);
			\draw[black,thick](2,0) edge[bend left= 60] node{\lmidarrow} (0,0);			
		\end{tikzpicture} \captionof{figure}{$C_2$}
			\end{minipage} & $H_{2}^{\infty +}$ &  [Theorem 4.6, \cite{Mandal}] \\
			\hline
			\begin{minipage}{0.3\linewidth}
				\centering
				\begin{tikzpicture}
				\draw[fill=black] (0,0) circle (2pt);
				\draw[fill=black] (2,0) circle (2pt);
				
				\begin{scope}[decoration={markings,
					mark=at position 0.25 with {\arrow[black,thick]{<}},
				}]
				\draw[black,postaction={decorate}, thick](0,0.5) circle (0.5);
				\end{scope}
				
				\begin{scope}[decoration={markings,
					mark=at position 0.25 with {\arrow[black,thick]{<}},
				}]
				\draw[black,postaction={decorate}, thick](2,0.5) circle (0.5);
				\end{scope}				
				\end{tikzpicture}
				\captionof{figure}{$L_1 \sqcup L_1$}
			\end{minipage}  & $H_2^{\infty +}$ & [Proposition 4.2, \cite{Mandalkms}]\\
			\hline
		\begin{minipage}{0.3\linewidth}
			\centering
			\begin{tikzpicture}   
			\draw[fill=black] (0,0) circle (2pt);
			\draw[fill=black] (2,0) circle (2pt);
			\draw[fill=black] (4,0) circle (2pt);

			\draw[black] (0,0)-- node{\lmidarrow}  (2,0);
			\draw[black] (2,0)-- node{\rmidarrow} (4,0);
			\end{tikzpicture}
			\captionof{figure}{$So_2$}
		\end{minipage}	 & $SH_2^{\infty +}$ &  [Proposition 3.4, \cite{max_permutational_symm}]  \\
			\hline
			\begin{minipage}{0.4\linewidth}
				\centering 
			\begin{tikzpicture}   
			\draw[fill=black] (0,0) circle (2pt);
			\draw[fill=black] (2,0) circle (2pt);
			\draw[fill=black] (3,0) circle (2pt);
			\draw[fill=black] (5,0) circle (2pt);
			
			\draw[black] (0,0)-- node{\rmidarrow}  (2,0);
			\draw[black] (3,0)-- node{\rmidarrow} (5,0);
			\end{tikzpicture}
				\captionof{figure}{$P_1 \sqcup P_1$}
			\end{minipage}  & $SH_2^{\infty +}$ & [Proposition 3.4, \cite{max_permutational_symm}] \\
			\hline
			 \begin{minipage}{0.3\linewidth} 
			 	\centering
			 	\begin{tikzpicture}   
			 	\draw[fill=black] (0,0) circle (2pt);
			 	\draw[fill=black] (2,0) circle (2pt);
			 	\draw[fill=black] (4,0) circle (2pt);
			 				 	
			 	\draw[black] (0,0)-- node{\rmidarrow}  (2,0);
			 	\draw[black] (2,0)-- node{\lmidarrow} (4,0);
			 	\end{tikzpicture} 
			 	\captionof{figure}{$Si_2$} \label{fig Si_2}
			  \end{minipage} & $U_2^{+}$ & [Proposition 3.3, \cite{max_permutational_symm}]\\ 
			\hline
		\end{tabular} \label{Table_U2+}
	\end{center}
Now, all the quantum symmetry groups $Q_{\tau}^{Lin}(\Gamma)$ listed in Table \ref{Table_U2+} (except $U_2^+$ itself) are not CQG isomorphic to $U_2^+$. To quickly justify this, we pass to the quotient of each of these CMQGs with their respective commutator ideal (say, $I_{com}$). Then $C(S^1)*C(S^1)/I_{com} \cong C(S^1 \times S^1)$, $H_2^{\infty +}/I_{com} \cong C(S^1 \wr S_2)$ and $SH_2^{\infty +}/I_{com} \cong C(S^1 \wr S_2)$. Clearly, the underlying compact topological groups $S^1 \times S^1$ and $S^1 \wr S_2$ are not homeomorphic to the unitary group $U(2)$. 
Hence, from Table \ref{Table_U2+}, one can conclude that the only possible choice of $\Gamma \in \mathcal{G}_{0,1}$ is $Si_2$, whose associated quantum symmetry group is isomorphic to $U_{2}^{+}$.  
\end{proof} 

\noindent Note: Recall that $Q_{\tau}^{Lin}(L_2) \cong U_2^+$ (by [Proposition 4.12, \cite{Mandal}]) and $Q_{\tau}^{Lin}(P_1^2) \cong U_2^+$ (by [Proposition 3.3,\cite{max_permutational_symm}]); however, neither $L_2$ nor $P_1^2$ belongs to $\mathcal{G}_{0,1}$. Here, $P_1^2$ denotes the graph with two edges having $A_{P_1^2}:=\begin{bmatrix}
    0 & 2\\
    0 & 0
\end{bmatrix}$.  \\

(2) However, for the CMQGs $\mystar\limits_{n} C(S^1)$ ($n \geq 1$), such a choice of graphs is not unique. Using Theorem IV.6 of \cite{rigidity}, one can ensure that there are several graphs with $n$ edges belonging to the family $\mathcal{G}$ (see Section III of \cite{rigidity}) having quantum symmetry $\mystar\limits_{n} C(S^1)$. For a concrete example, observe that for non-isomorphic graphs, $P_{n}$ and $T_{n}$ satisfy $Q_{\tau}^{Lin}(P_{n}) \cong \mystar\limits_{n} C(S^1) \cong Q_{\tau}^{Lin}(T_{n})$. \\

\begin{figure}[htpb]
	\begin{minipage}{0.4\linewidth}
		\centering
		\begin{tikzpicture}[scale=0.8]
		\draw[fill=black] (0,0) circle (2pt) node[anchor=north]{$v_{1}$};
		\draw[fill=black] (2,0) circle (2pt) node[anchor=north]{$v_{2}$};
		\draw[fill=black] (4,0) circle (2pt) node[anchor=north]{$v_{3}$};
		\draw[fill=black] (6,0) circle (2pt) node[anchor=north]{$v_{n}$};
		\draw[fill=black] (8,0) circle (2pt) node[anchor=north]{$v_{n+1}$};
		
		\draw[black,thick] (0,0)-- node{\rmidarrow}  node[above]{$e_{12}$} (2,0);
		\draw[black,thick] (2,0)-- node{\rmidarrow} node[above]{$e_{23}$} (4,0);
		\draw[loosely dotted, black,thick] (4.5,0)-- (5.5,0);
		\draw[black,thick] (6,0)--node{\rmidarrow}  node[above]{$e_{n(n+1)}$} (8,0);
		\end{tikzpicture}
		\captionof{figure}{$P_n$}
	\end{minipage}
\hspace{1.5cm}
\begin{minipage}{0.4\linewidth}
	\centering
	\begin{tikzpicture}[scale=0.8]
	\draw[fill=black] (0,0) circle (2pt) node[anchor=north]{$v_{1}$};
	\draw[fill=black] (2,0) circle (2pt) node[anchor=north]{$v_{2}$};
	\draw[fill=black] (4,0) circle (2pt) node[anchor=north]{$v_{3}$};
	\draw[fill=black] (6,0) circle (2pt) node[anchor=north]{$v_{n-1}$};
	\draw[fill=black] (8,0) circle (2pt) node[anchor=north]{$v_{n}$};
	
	\begin{scope}[decoration={markings,
		mark=at position 0.25 with {\arrow[black,thick]{<}},
	}]
	\draw[black,postaction={decorate}, thick]  (0,0.5) circle (0.5) node[above=0.5]{$e_{11}$};
	\end{scope}
	
	\draw[black,thick] (0,0)-- node{\rmidarrow}  node[above]{$e_{12}$} (2,0);
	\draw[black,thick] (2,0)-- node{\rmidarrow} node[above]{$e_{23}$} (4,0);
	\draw[loosely dotted, black,thick] (4.5,0)-- (5.5,0);
	\draw[black,thick] (6,0)--node{\rmidarrow}  node[above]{$e_{(n-1)n}$} (8,0);
	\end{tikzpicture}
	\captionof{figure}{$T_{n}$}
\end{minipage}
\end{figure}

These observations lead to the following general question:\\ 
Can one identify all such quantum symmetry groups $Q \in \mathcal{Q}_{JM}$ for which there exists a unique graph $\Gamma_0$ having quantum symmetry $Q_{\tau}^{Lin}(\Gamma_0) \cong Q$?\\
Clearly, such a family of quantum symmetry groups is nonempty by Lemma \ref{prop U_2^+}.

\subsection{Complement graphs:} \label{sec_comp. graph}
Since we consider directed graphs without multiple edges (but possibly with loops), we define the {\bf complement graph} $\Gamma^{c}$ of a graph $\Gamma$ to be the directed graph whose adjacency matrix $A_{\Gamma^c}$ is given by $J_n - A_{\Gamma}$. \\

\noindent $\blacksquare$ Now, if $u$ denotes the canonical fundamental representation of the quantum symmetry group $QAut_{Ban}(\Gamma)$, then clearly $uA_{\Gamma}=A_{\Gamma}u \iff uA_{\Gamma^c}=A_{\Gamma^c}u $. Thus, universal properties of respective quantum symmetry groups ensure that $QAut_{Ban}(\Gamma) \cong QAut_{Ban}(\Gamma^c)$. \\

\noindent $\blacktriangleright$ In contrast, the analogous isomorphism does not hold in general in the graph $C^*$-algebraic setting. That is, $Q_{\tau}^{Lin}(\Gamma)$ may not be isomorphic to $Q_{\tau}^{Lin}(\Gamma^c)$ in general. For instance, take $\Gamma$ to be $P_1$, then $\Gamma^c \cong \bar{L}_3$, where $A_{\bar{L}_3}:=\begin{bmatrix}
    1 & 1\\
    0 & 1
\end{bmatrix}$ (see [FIG. 4, \cite{rigidity}]). Therefore, by [Theorem IV.6, \cite{rigidity}], $Q_{\tau}^{Lin}(P_1) \cong C(S^1) \ncong \mystar\limits_{3} C(S^1) \cong Q_{\tau}^{Lin}(\bar{L}_3)$.\\
However, when making such a comparison, we need to be careful that both $\Gamma$ and $\Gamma^c$ belong to $\mathcal{G}_{0,1}$ (to be precise, both $\Gamma$ and $\Gamma^c$ must be graphs without isolated vertices so that the notion of quantum symmetry groups $Q_{\tau}^{Lin}(\Gamma)$ and $Q_{\tau}^{Lin}(\Gamma^c)$ can be given).\\
Interestingly, using the family $\mathcal{G}$ described in [Section III, \cite{rigidity}], one can construct an infinite family of graphs $\Gamma_n$ (resp. $\Lambda_n$) $\in \mathcal{G}_{0,1}$ satisfying $Q_{\tau}^{Lin}(\Gamma_n) \ncong Q_{\tau}^{Lin}(\Gamma_n^c)$ (resp. $Q_{\tau}^{Lin}(\Lambda_n) \cong Q_{\tau}^{Lin}(\Lambda_n^c)$).\\

(1) Consider the graphs $\Gamma_n := L_{2n-1} \in \mathcal{G}$ ($n >1$) whose adjacency matrix $A_{L_{2n-1}}$ is given by $(A_{L_{2n-1}})_{ij}:= \begin{cases}
1 & if ~ i \leq j \\
0 & if ~ i>j
\end{cases}$  (see [FIG. 3, \cite{rigidity}]). Then the complement graph $\Gamma^c := L_{2n-1}^c$ also belongs to $\mathcal{G}$, since its adjacency matrix $(A_{L_{2n-1}^c})_{ij} := \begin{cases}
0 & if ~ i \leq j \\
1 & if ~ i>j
\end{cases}$. By [Theorem IV.6, \cite{rigidity}], we obtain $$Q_{\tau}^{Lin}(L_{2n-1}) \cong \mystar\limits_{|E(\Gamma_n)|} C(S^1) \ncong \mystar\limits_{|E(\Gamma_n^c)|} C(S^1) \cong Q_{\tau}^{Lin}(L_{2n-1}^c),$$ as $|E(\Gamma_n)| - |E(\Gamma_n^c)|=n$. \\

(2) Consider the graphs $\Lambda_n \in \mathcal{G}$ with $n$ vertices (where $n > 1$ is multiple of 4) whose adjacency matrix $A_{\Lambda_n}$ is given by $$(A_{\Lambda_n})_{ij}:= \begin{cases}
1 & if ~ i < j \\
1 & if ~ i=j \in \{1, 2,..., n/4 \} \cup \{3n/4+1,..., n-1, n \}\\
0 & otherwise \\ 
\end{cases}.$$ 
Clearly, by constriction of $\Lambda_n$, its complement graph $\Lambda_n^c \in \mathcal{G}$ is not isomorphic to $\Lambda_n$; however, $|E(\Lambda_n)| = |E(\Lambda_n)|$ (see Figures \ref{Lambda_8} \& \ref{Lambda_8^c}). Thus, again by [Theorem IV.6, \cite{rigidity}], we conclude that $$Q_{\tau}^{Lin}(\Lambda_n) \cong \mystar\limits_{|E(\Lambda_n)|} C(S^1) \cong Q_{\tau}^{Lin}(\Lambda_n^c). $$  

 \begin{center}
	\begin{figure}[htb]   
		\begin{tikzpicture}
		\draw[fill=black] (0,0) circle (2pt) node[anchor=south]{$v_{1}$};
		\draw[fill=black] (2,0) circle (2pt) node[anchor=south]{$v_{2}$};
		\draw[fill=black] (4,0) circle (2pt) node[anchor=south]{$v_{3}$};
		\draw[fill=black] (6,0) circle (2pt) node[anchor=south]{$v_{4}$};
		\draw[fill=black] (8,0) circle (2pt) node[anchor=south]{$v_{5}$};
		\draw[fill=black] (10,0) circle (2pt) node[anchor=south]{$v_{6}$};
		\draw[fill=black] (12,0) circle (2pt) node[anchor=south]{$v_{7}$};
		\draw[fill=black] (14,0) circle (2pt) node[anchor=south]{$v_{8}$};

		\draw[black,thick](0,0.5) circle (0.5) node[above=0.25]{\rmidarrow};
		\draw[black,thick](2,0.5) circle (0.5) node[above=0.25]{\rmidarrow};
		\draw[black,thick](12,0.5) circle (0.5)  node[above=0.25]{\rmidarrow};
		\draw[black,thick](14,0.5) circle (0.5)  node[above=0.25]{\rmidarrow};

		\draw[black,thick] (0,0)-- node{\rmidarrow}  (2,0);
		\draw[black,thick] (2,0)-- node{\rmidarrow}  (4,0);
		\draw[black,thick] (4,0)--node{\rmidarrow}   (6,0);
		\draw[black,thick] (6,0)--node{\rmidarrow}   (8,0);
		\draw[black,thick] (8,0)--node{\rmidarrow}   (10,0);
		\draw[black,thick] (10,0)--node{\rmidarrow}  (12,0);
		\draw[black,thick] (12,0)--node{\rmidarrow}  (14,0);
		
		\draw[black,thick] (0,0) edge[bend left= -25] node{\rmidarrow}  (4,0);
		\draw[black,thick] (0,0) edge[bend left= -35] node{\rmidarrow}  (6,0);
		\draw[black,thick] (0,0) edge[bend left= -45] node{\rmidarrow}  (8,0);
		\draw[black,thick] (0,0) edge[bend left= -55] node{\rmidarrow}  (10,0);
		\draw[black,thick] (0,0) edge[bend left= -65] node{\rmidarrow}  (12,0);
		\draw[black,thick] (0,0) edge[bend left= -75] node{\rmidarrow}  (14,0);

		\draw[black,thick] (2,0) edge[bend left= -25] node{\rmidarrow}  (6,0);
		\draw[black,thick] (2,0) edge[bend left= -35] node{\rmidarrow}  (8,0);
		\draw[black,thick] (2,0) edge[bend left= -45] node{\rmidarrow}  (10,0);
		\draw[black,thick] (2,0) edge[bend left= -55] node{\rmidarrow}  (12,0);
		\draw[black,thick] (2,0) edge[bend left= -65] node{\rmidarrow}  (14,0);
		
		\draw[black,thick] (4,0) edge[bend left= -25] node{\rmidarrow}  (8,0);
		\draw[black,thick] (4,0) edge[bend left= -35] node{\rmidarrow}  (10,0);
		\draw[black,thick] (4,0) edge[bend left= -45] node{\rmidarrow}  (12,0);
		\draw[black,thick] (4,0) edge[bend left= -55] node{\rmidarrow}  (14,0);
		
       \draw[black,thick] (6,0) edge[bend left= -25] node{\rmidarrow}  (10,0);
        \draw[black,thick] (6,0) edge[bend left= -35] node{\rmidarrow} (12,0);
        \draw[black,thick] (6,0) edge[bend left= -45] node{\rmidarrow} (14,0);
        
       \draw[black,thick] (8,0) edge[bend left= -25] node{\rmidarrow}   (12,0);
       \draw[black,thick] (8,0) edge[bend left= -35] node{\rmidarrow}   (14,0);
   
       \draw[black,thick] (10,0) edge[bend left= -25] node{\rmidarrow}  (14,0);             				
		\end{tikzpicture}
		\caption{$\Lambda_8$}    \label{Lambda_8}
	\end{figure}
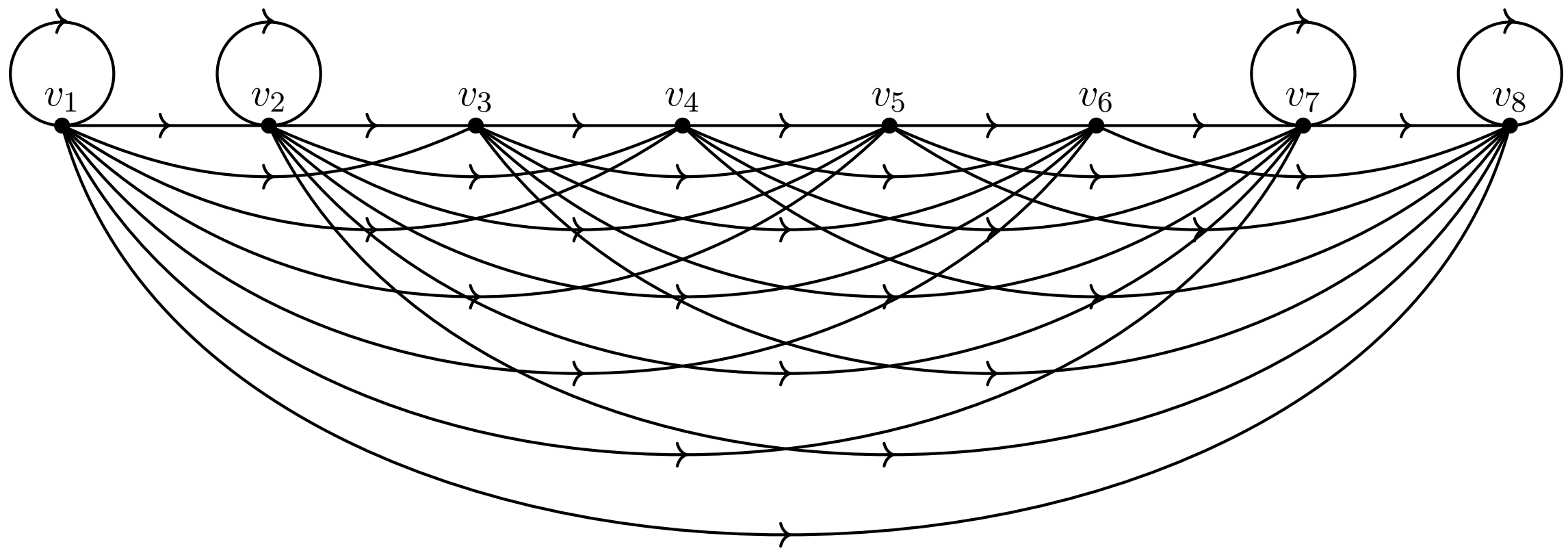
\end{center}

\begin{center}
	\begin{figure}[htb]   
		\begin{tikzpicture}
		\draw[fill=black] (0,0) circle (2pt) node[anchor=south]{$v_{1}$};
		\draw[fill=black] (2,0) circle (2pt) node[anchor=south]{$v_{2}$};
		\draw[fill=black] (4,0) circle (2pt) node[anchor=south]{$v_{3}$};
		\draw[fill=black] (6,0) circle (2pt) node[anchor=south]{$v_{4}$};
		\draw[fill=black] (8,0) circle (2pt) node[anchor=south]{$v_{5}$};
		\draw[fill=black] (10,0) circle (2pt) node[anchor=south]{$v_{6}$};
		\draw[fill=black] (12,0) circle (2pt) node[anchor=south]{$v_{7}$};
		\draw[fill=black] (14,0) circle (2pt) node[anchor=south]{$v_{8}$};
				
		\draw[black,thick](4,0.5) circle (0.5) node[above=0.25]{\lmidarrow};
		\draw[black,thick](6,0.5) circle (0.5)  node[above=0.25]{\lmidarrow};
		\draw[black,thick](8,0.5) circle (0.5)  node[above=0.25]{\lmidarrow};
		\draw[black,thick](10,0.5) circle (0.5)  node[above=0.25]{\lmidarrow};

		\draw[black,thick] (0,0)-- node{\lmidarrow}  (2,0);
		\draw[black,thick] (2,0)-- node{\lmidarrow}  (4,0);
		\draw[black,thick] (4,0)--node{\lmidarrow}   (6,0);
		\draw[black,thick] (6,0)--node{\lmidarrow}   (8,0);
		\draw[black,thick] (8,0)--node{\lmidarrow}   (10,0);
		\draw[black,thick] (10,0)--node{\lmidarrow}  (12,0);
		\draw[black,thick] (12,0)--node{\lmidarrow}  (14,0);
		
		\draw[black,thick] (0,0) edge[bend left= -25] node{\lmidarrow}  (4,0);
		\draw[black,thick] (0,0) edge[bend left= -35] node{\lmidarrow}  (6,0);
		\draw[black,thick] (0,0) edge[bend left= -45] node{\lmidarrow}  (8,0);
		\draw[black,thick] (0,0) edge[bend left= -55] node{\lmidarrow}  (10,0);
		\draw[black,thick] (0,0) edge[bend left= -65] node{\lmidarrow}  (12,0);
		\draw[black,thick] (0,0) edge[bend left= -75] node{\lmidarrow}  (14,0);
		
		\draw[black,thick] (2,0) edge[bend left= -25] node{\lmidarrow}  (6,0);
		\draw[black,thick] (2,0) edge[bend left= -35] node{\lmidarrow}  (8,0);
		\draw[black,thick] (2,0) edge[bend left= -45] node{\lmidarrow}  (10,0);
		\draw[black,thick] (2,0) edge[bend left= -55] node{\lmidarrow}  (12,0);
		\draw[black,thick] (2,0) edge[bend left= -65] node{\lmidarrow}  (14,0);
		
		\draw[black,thick] (4,0) edge[bend left= -25] node{\lmidarrow}  (8,0);
		\draw[black,thick] (4,0) edge[bend left= -35] node{\lmidarrow}  (10,0);
		\draw[black,thick] (4,0) edge[bend left= -45] node{\lmidarrow}  (12,0);
		\draw[black,thick] (4,0) edge[bend left= -55] node{\lmidarrow}  (14,0);
		
		\draw[black,thick] (6,0) edge[bend left= -25] node{\lmidarrow}  (10,0);
		\draw[black,thick] (6,0) edge[bend left= -35] node{\lmidarrow} (12,0);
		\draw[black,thick] (6,0) edge[bend left= -45] node{\lmidarrow} (14,0);
		
		\draw[black,thick] (8,0) edge[bend left= -25] node{\lmidarrow}   (12,0);
		\draw[black,thick] (8,0) edge[bend left= -35] node{\lmidarrow}   (14,0);
		
		\draw[black,thick] (10,0) edge[bend left= -25] node{\lmidarrow}  (14,0);	
		\end{tikzpicture}
		\caption{$\Lambda_8^c$}   \label{Lambda_8^c}
	\end{figure}
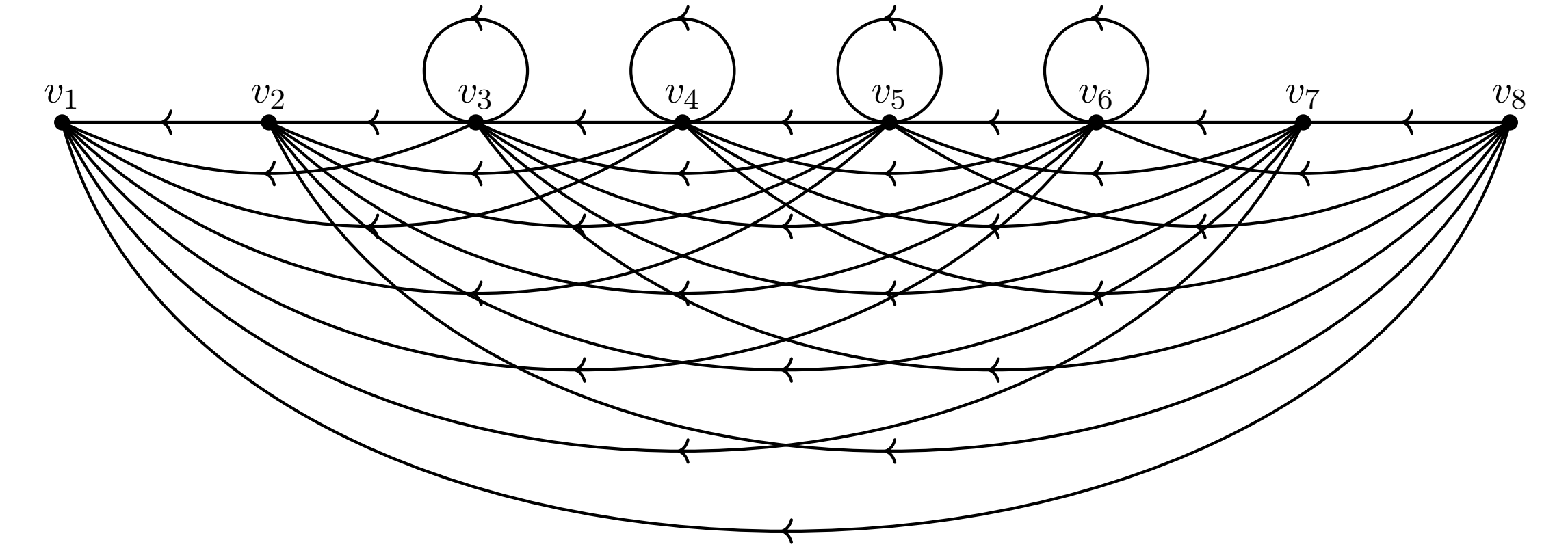
\end{center}

\begin{rem}
It is also possible to find a graph $\Gamma \in \mathcal{G}_{0,1}$ with a non-minimal quantum symmetry group $Q \in \mathcal{Q}_{JM}$ (i.e., $Q$ is not of the form of free product $C(S^1)$) such that $Q_{\tau}^{Lin}(\Gamma) \cong Q_{\tau}^{Lin}(\Gamma^c)$. For example, let $\Gamma = L_1 \sqcup L_1$. Then $\Gamma^c = C_2$. Hence, using [Theorem 4.6, \cite{Mandal}] and [Proposition 4.2, \cite{Mandalkms}], one can conclude that $Q_{\tau}^{Lin}(L_1 \sqcup L_1) \cong H_2^{\infty +} \cong Q_{\tau}^{Lin}(C_2)$. Now, this naturally leads to the following question:\\
Can one find an infinite family of graphs $\Gamma$ such that $Q_{\tau}^{Lin}(\Gamma) \cong Q_{\tau}^{Lin}(\Gamma^c) \ncong \mystar\limits_{|E(\Gamma)|} C(S^1)$? \\  
Corollary \ref{edge equal} at least ensures that any such graphs $\Gamma$ satisfy $|E(\Gamma)|=|E(\Gamma^c)|=\frac{|V(\Gamma)|^2}{2}$ (since we allow loops in a graph). \\ 
\end{rem}

\subsection{Disjoint union of non-isomorphic graphs} \label{sec_free prod}
We compare the behaviour of two notions of quantum symmetry groups of disjoint union of non-isomorphic graphs in $\mathcal{G}_{0,1}$ by investigating whether the corresponding quantum symmetry decomposes into the free product of individual quantum symmetries. \\[0.1cm]
If $\Gamma_1$ and $\Gamma_2 \in \mathcal{G}_{0,1}$ are weakly connected quantum non-isomorphic graphs, then $QAut_{Ban}(\Gamma_1 \sqcup \Gamma_2) \cong QAut_{Ban}(\Gamma_1) * QAut_{Ban}(\Gamma_2)$ [Proposition \ref{graph_free prod}]. But the similar isomorphism does not hold in general in the context of graph $C^*$-algebra when $QAut_{Ban}(\cdot)$ is replaced by $Q_{\tau}^{Lin}(\cdot)$ [Proposition \ref{graph C^*alg_free prod not hold}].\\[0.2cm]
Consequently, it holds that given $Q \in \mathcal{Q}_{Ban}$, there always exist infinitely many mutually non-isomorphic weakly connected graphs $\{\Gamma_i\}_{i \in \mathbb{N}} \in \mathcal{G}_{0,1}$ [Theorem \ref{graph_infinitely many}] (as well as disconnected graphs $\{\Gamma^{dis}_i\}_{i \in \mathbb{N}} \in \mathcal{G}_{0,1}$ [Corollary \ref{discon. inf. many graphs}]) such that $QAut_{Ban}(\Gamma_i) \cong Q \cong QAut_{Ban}(\Gamma^{dis}_i)$ for all $i \in \mathbb{N}$. In contrast, in the graph $C^*$-algebraic scenario, for the quantum symmetry group $U_2^+ \in \mathcal{Q}_{JM}$, there is no disconnected graph $\Gamma$ having quantum symmetry $Q_{\tau}^{Lin}(\Gamma) \cong U_2^+$ (by Proposition \ref{prop U_2^+}).  \\[0.3cm]
$\blacksquare$ In the graph scenario, we have the following:
\begin{prop} \label{graph_free prod}
If $\Gamma_1$ and $\Gamma_2$ are weakly connected, quantum non-isomorphic directed graphs with $A_{\Gamma_i} \in M_n(\{0,1\})$, then $QAut_{Ban}(\Gamma_1 \sqcup \Gamma_2) \cong QAut_{Ban}(\Gamma_1) * QAut_{Ban}(\Gamma_2)$
\end{prop}
An analogous result has already been established for the simple undirected graphs in [Lemma 6.4, \cite{Qauttrees}] (also refer to [Theorem 4.5, \cite{Lupini}] and [Corollary 7.1.4, \cite{SSthesis}]). We will just ensure that the similar proof carries over in our setting. We briefly outline the proof below, mimicked from [Lemmas 6.3 \& 6.4, \cite{Qauttrees}].\\[0.1cm]
The first step is to show that under the given hypothesis, $u_{st}=0$ for all $s \in V(\Gamma_1)$ and $t \in V(\Gamma_2)$, where $u$ denotes the canonical fundamental representation of $QAut_{Ban}(\Gamma_1 \sqcup \Gamma_2)$.
To this end, the key point is to verify that the analogue of [Lemma 6.2, \cite{Qauttrees}] remains valid in our setting; that is, $u_{xy} u_{x'y'}=0$ if both $x,x'$ belong to the same connected component, but $y$ and $y'$ belong to different components. Since each $\Gamma_i$ is weakly connected, then $x =: x_0 \sim x_1 \sim x_2 \sim \cdots \sim x_n:=x'$ for some vertices $x_i$ in that fixed component graph. Now,
\begin{align*}
u_{xy}u_{x'y'} &= u_{xy} \left(\sum\limits_{y_i \in V(\Gamma_1 \sqcup \Gamma_2)} u_{x_1 y_1} \right) \left(\sum\limits_{y_2 \in V(\Gamma_1 \sqcup \Gamma_2)} u_{x_2 y_2} \right) \cdots \left(\sum\limits_{y_{n-1} \in V(\Gamma_1 \sqcup \Gamma_2)} u_{x_{n-1} y_{n-1}} \right)u_{x'y'}\\
& = \sum\limits_{y_1,...y_{n-1} \in V(\Gamma_1 \sqcup \Gamma_2)} u_{xy}  u_{x_1 y_1} \cdots  u_{x_{n-1} y_{n-1}} u_{x'y'}
\end{align*}
Since $y_0:= y$ and $y_n:= y'$ are in different components, for some $i \in \{0\} \sqcup [n]$, neither $y_i \rightarrow y_{i+1}$ nor $y_{i+1} \rightarrow y_{i}$. Thus, by [Lemma 3.1.1, \cite{Fulton}], if $x_{i} \rightarrow x_{i+1}$, then directly $u_{{x_i}{y_i}}u_{x_{i+1}y_{i+1}}=0$; and if $x_{i+1} \rightarrow x_{i}$, then $u_{x_{i+1}y_{i+1}} u_{{x_i}{y_i}}=0 \implies u_{{x_i}{y_i}} u_{x_{i+1}y_{i+1}} = (u_{x_{i+1}y_{i+1}} u_{{x_i}{y_i}})^*=0$. Hence, $u_{xy}u_{x'y'}=0$. 
The rest of the proof follows exactly the same lines as that of Lemma 6.3 of \cite{Qauttrees}, since it consists purely of $C^*$-algebraic calculations and uses only the fact that $u$ is a magic unitary matrix. \\[0.1cm]
The second step (analogue of [Lemma 6.4, \cite{Qauttrees}]) shows that $QAut_{Ban}(\Gamma_1 \sqcup \Gamma_2) \cong QAut_{Ban}(\Gamma_1) * QAut_{Ban}(\Gamma_2)$, where the CQG isomorphism sends $u$ to the block diagonal matrix whose diagonal blocks are the fundamental representations of $QAut_{Ban}(\Gamma_1)$ and $QAut_{Ban}(\Gamma_2)$ respectively. Since $u$ is decomposed into diagonal blocks 
$\begin{bmatrix}
u^{(1)} & 0 \\
0 &  u^{(2)}
\end{bmatrix}$ and $A_{\Gamma}:= \begin{bmatrix}
A_{\Gamma_1} & 0 \\
0 & A_{\Gamma_2}
\end{bmatrix}$, $A_{\Gamma}u= uA_{\Gamma}$ is equivalent to the equations $A_{\Gamma_i}u^{(i)}=u^{(i)}A_{\Gamma_i}$ (for $i=1,2$). The results therefore follow from universal properties of the corresponding quantum symmetry groups.\\[0.2cm]

In Theorem \ref{graph_infinitely many}, we have seen that given a directed graph $\Gamma$ (without multiple edges) on $n (>1)$ vertices, there always exist infinitely many mutually non-isomorphic weakly connected directed graphs having their quantum automorphism group invariant. As a consequence of Proposition \ref{graph_free prod}, we can ensure the same for the disconnected graph as well.
\begin{cor} \label{discon. inf. many graphs}
    Given a directed graph $\Gamma$ (without multiple edges) on $n (>1)$ vertices, there always exists an infinite family of mutually non-isomorphic disconnected directed graphs  $\{\Gamma^{dis}_i\}$ such that $QAut_{Ban}(\Gamma^{dis}_i) \cong QAut_{Ban}(\Gamma)$ for all $i$.
\end{cor}
\begin{proof}
   It is easy to verify that $QAut_{Ban}(P_{n}) \cong \mathbb{C}$ for all $n \in \mathbb{N}$. Let $\Gamma$ be weakly connected. Choose an integer $m > |V(\Gamma)|$. Consider the disconnected graphs $\Gamma^{dis}_{i}:= \Gamma \sqcup P_i$ for $i > m$. Since $\Gamma$ is quantum non-isomorphic to $P_i$, Proposition \ref{graph_free prod} implies that $QAut_{Ban}(\Gamma^{dis}_{i}) \cong QAut_{Ban}(\Gamma)*\mathbb{C} \cong QAut_{Ban}(\Gamma)$ for all $i>m$. Now, if $\Gamma$ itself is disconnected, then directly applying the constructions mentioned in STEP-2 of Theorem \ref{graph_infinitely many} to $\Gamma$, one gets infinitely many $\{\Gamma^{dis}_i\}_{i>m}$ having $QAut_{Ban}(\Gamma^{dis}_i) \cong QAut_{Ban}(\Gamma)$ for all $i$. Indeed, the graph constructed via STEP-2 of Theorem \ref{graph_infinitely many} remains disconnected if the original graph is disconnected.
\end{proof}

\noindent $\blacktriangleright$ In the graph $C^*$-algebraic setting, it can be shown that if $\Gamma_1$ and $\Gamma_2$ are two directed graphs without isolated vertices, then CMQG $(Q_{\tau}^{Lin}(\Gamma_1) * Q_{\tau}^{Lin}(\Gamma_2), q^{(1)} \oplus q^{(2)})$ acts on $C^*(\Gamma)$ linearly, faithfully and $\tau$-preserving manner via the action $\alpha: C^*(\Gamma_1 \sqcup \Gamma_2) \to C^*(\Gamma_1 \sqcup \Gamma_2) \otimes (Q_{\tau}^{Lin}(\Gamma_1) * Q_{\tau}^{Lin}(\Gamma_2))$ defined by
\begin{align*}
\alpha(S_e) &= \sum\limits_{f \in E(\Gamma_1)} S_f \otimes q^{(1)}_{fe} \text{ for all } e \in E(\Gamma_1) \\
\text{ \ \ \ and \ \ \ } \alpha(S_e) &= \sum\limits_{f \in E(\Gamma_2)} S_f \otimes q^{(2)}_{fe} \text{ for all } e \in E(\Gamma_2).
\end{align*} 
Therefore, $Q_{\tau}^{Lin}(\Gamma_1) * Q_{\tau}^{Lin}(\Gamma_2)$ is always a quantum subgroup of $Q_{\tau}^{Lin}(\Gamma_1 \sqcup \Gamma_2)$.\\[0.2cm]
However, as mentioned earlier in the counter-example of [Section 3.1, \cite{cuntz_KM}], $Q_{\tau}^{Lin}(\Gamma_1 \sqcup \Gamma_2) \ncong Q_{\tau}^{Lin}(\Gamma_1) * Q_{\tau}^{Lin}(\Gamma_2)$ in general, even for quantum non-isomorphic graphs. The following proposition guarantees that one can actually find an infinite pair of weakly connected quantum non-isomorphic graphs for which the above isomorphism does not hold. 
\begin{prop} \label{graph C^*alg_free prod not hold}
	If $\Gamma_1, \Gamma_2$ are two graphs in the family $\mathcal{I}_1$, where $$ \mathcal{I}_1 :=\{\Gamma : \Gamma \text{ does not contain any loop and for each } e \in E(\Gamma), r(e) \text{ is a sink of indegree 1}\}  \subset \mathcal{G}_{0,1},$$ then $Q_{\tau}^{Lin}(\Gamma_1 \sqcup \Gamma_2) \ncong Q_{\tau}^{Lin}(\Gamma_1) * Q_{\tau}^{Lin}(\Gamma_2)$.
\end{prop}
\begin{proof}
	Observe that if $\Gamma_1, \Gamma_2 \in \mathcal{I}_{1}$, then $\Gamma_1 \sqcup \Gamma_2 \in \mathcal{I}_{1}$. Therefore, by [Proposition 3.4, \cite{max_permutational_symm}], we have $Q_{\tau}^{Lin}(\Gamma_1 \sqcup \Gamma_2) \cong SH_{m+n}^{\infty +}$, $Q_{\tau}^{Lin}(\Gamma_1) \cong SH_{m}^{\infty +}$, and $Q_{\tau}^{Lin}(\Gamma_2) \cong SH_{n}^{\infty +}$, where $m:=|E(\Gamma_1)|$ and $n:=|E(\Gamma_2)|$. Now, we claim that $SH_{m+n}^{\infty +} \ncong SH_{m}^{\infty +} * SH_{n}^{\infty +}$. Suppose for the sake of contradiction that such an isomorphism exists. Then passing to the quotient by their respective commutator ideals, we obtain $C(S^1 \wr S_{m+n}) \cong C(S^1 \wr S_m) \otimes C(S^1 \wr S_n) \cong C((S^1 \wr S_m) \times (S^1 \wr S_n))$. Thus, comparing the underlying compact topological spaces associated with these commutative $C^*$-algebras, it follows that $(S^1)^{m+n} \times S_{m+n}$ is homeomorphic $ ((S^1)^{m} \times S_{m}) \times ((S^1)^{n} \times S_{n}) $, which is impossible. Indeed, the first one has $(m+n)!$ many connected components, whereas the latter has $m!n! $ many connected components.
\end{proof}        
\begin{figure}
	\begin{minipage}[c]{0.4\linewidth}
		\centering
		\begin{tikzpicture}
		\draw[fill=black] (0,0) circle (2pt) node[anchor=north]{$v$};
		\draw[fill=black] (2,0) circle (2pt) node[anchor=north]{$v_{n}$};
		\draw[fill=black] (1.414,1.414) circle (2pt) node[anchor=north]{$v_{1}$};
		\draw[fill=black] (0,2) circle (2pt) node[anchor=south]{$v_{2}$};
		\draw[fill=black] (-1.414,1.414) circle (2pt) node[anchor=north]{$v_{3}$};
		
		\draw[black] (0,0)-- node{\rmidarrow}  node[below]{$e_{n}$} (2,0);
		\draw[black] (0,0)-- node{\Nwmidarrow}  node[right]{$e_{1}$} (1.414,1.414);
		\draw[black] (0,0)-- node{\Nmidarrow}  node[right]{$e_{2}$} (0,2);
		\draw[black] (0,0)-- node{\Nemidarrow}  node[above]{$e_{3}$} (-1.414,1.414);
		
		\draw [loosely dotted,thick,domain=150:340] plot ({1.5*cos(\x)}, {1.5*sin(\x)});
		\end{tikzpicture}
		\caption{$So_n$} \label{Son}
	\end{minipage}
\end{figure}
\begin{rem}
	Note that any graph in the family $\mathcal{I}_{1}$ is weakly connected iff it is isomorphic to $So_n$ (see Figure \ref{Son}), where $n$ denotes the number of edges. Moreover, for any $m,n \in \mathbb{N}$ with $m \neq n$, the graphs $So_m$ and $So_n$ are quantum non-isomorphic, while $Q_{\tau}^{Lin}(So_m \sqcup So_n) \ncong Q_{\tau}^{Lin}(So_m) * Q_{\tau}^{Lin}(So_n)$. 
\end{rem}

On the other hand, it is also possible to produce a family of graphs from $\mathcal{G}$ (as described in [Section III, \cite{rigidity}]), namely, the path graphs $P_n$, whose associated graph $C^*$-algebras $C^*(P_n)$ are isomorphic to the matrix algebras $M_{n+1}(\mathbb{C})$. These graphs satisfy the relation $Q_{\tau}^{Lin}(P_n \sqcup P_m) \cong Q_{\tau}^{Lin}(P_n) * Q_{\tau}^{Lin}(P_m)$ for $m \neq n$. More generally, we have the following.
\begin{prop} \label{matix free prod}
    	Let $\{M_{n_i +1}(\mathbb{C})\}_{i=1}^{m}$ be a finite family of matrix algebras where all $n_i$'s are distinct. Then $ Q_{\tau}^{Lin}(\sqcup_{i=1}^{m} ~ P_{n_i}) \cong \mystar\limits_{n} C(S^1) \cong \mystar\limits_{i=1}^{m}Q_{\tau}^{Lin}(P_{n_i}) $, where $n=\sum\limits_{i=1}^{m} n_i$.
\end{prop}
\noindent We omit the proof here and refer to Theorem 5.1.4 in Chapter 5 of \cite{Ujjal_thesis}.

\subsection{Behavior of minimal quantum symmetry groups} \label{sec_behav. min. q. symm.}
Recall from [Proposition 3.11, \cite{rigidity}] that $\left(\mystar\limits_{|E(\Gamma)|} C(S^1), u \right)$ always acts linearly, faithfully and $\tau$-preserving manner on a graph $C^*$-algebra $C^*(\Gamma)$, where $u:= diag(u_1, u_2, \cdots, u_{|E(\Gamma)|})$ is the canonical fundamental representation with diagonal entries unitary elements $u_i$'s.  
In view of that, we say the {\bf quantum symmetry of $C^*(\Gamma)$ is minimal} when the canonical fundamental representation associated to its quantum symmetry group $(Q_{\tau}^{Lin}(\Gamma),q)$ takes the form $q:=diag(q_1, q_2, \cdots, q_{|E(\Gamma)|})$ where $q_i$'s are unitary elements. Obviously, in that case, $(Q_{\tau}^{Lin}(\Gamma),q)$ is CQG isomorphic to $\left(\mystar\limits_{|E(\Gamma)|} C(S^1), u \right)$.\\

\noindent In contrast, we can say that the quantum symmetry of $\Gamma$ is minimal whenever $\Gamma$ is {\bf quantum asymmetric}; that is, the canonical fundamental representation $u$ associated to its quantum symmetry group $(QAut_{Ban}(\Gamma),u)$ is equal to $I_{|V(\Gamma)|}$ or, equivalently, $(QAut_{Ban}(\Gamma),u) \cong \mathbb{C}$.\\

From this perspective, we compare the two notions of quantum symmetry. First, if two graphs from $\mathcal{G}_{0,1}$ have a distinct number of edges but the quantum symmetries of their respective graph $C^*$-algebras are both minimal, then also their quantum symmetry groups are necessarily distinct (since $\mystar\limits_{m} C(S^1) \ncong  \mystar\limits_{n} C(S^1)$ for $m \neq n$). On the other hand, if any graphs from $\mathcal{G}_{0,1}$ have `minimal quantum symmetry in a graph-theoretic sense' (formally, they are quantum asymmetric), then their quantum symmetry groups are both isomorphic to $\mathbb{C}$.\\[0.2cm]
 More interestingly, for a graph $\Gamma \in \mathcal{G}_{0,1}$, the following relation holds between these two notions of quantum symmetry, namely $Q_{\tau}^{Lin}(\Gamma)$ and  $QAut_{Ban}(\Gamma)$.
\begin{enumerate}
	\item If quantum symmetry of $C^*(\Gamma)$ is minimal, then $\Gamma$ is always quantum asymmetric.
	\item However, being $\Gamma$ quantum asymmetric does not necessarily imply that the quantum symmetry of $C^*(\Gamma)$ is minimal.
\end{enumerate}

Let $q$ and $u$ denote the canonical fundamental representations of $Q_{\tau}^{Lin}(\Gamma)$ and $QAut_{Ban}(\Gamma)$ respectively. With this convention, we show (1) in the following proposition.
\begin{prop} \label{rigid imply q.assym}
	If quantum symmetry of $C^*(\Gamma)$ is minimal (i.e., $q$ is diagonal with unitary diagonal entries), then $\Gamma$ is always quantum asymmetric (i.e., $u=I_{|V(\Gamma)|} \in M_{|V(\Gamma)|}(QAut_{Ban}(\Gamma))$).
\end{prop}
\begin{proof}
	It is known from [Theorem 3.10, \cite{Mandal}] that $(QAut_{Ban}(\Gamma),u)$ acts on $C^*(\Gamma)$ via the linear faithful action $\alpha: C^*(\Gamma) \to C^*(\Gamma) \otimes QAut_{Ban}(\Gamma)$ defined by $$\alpha(S_e)= \sum\limits_{f \in E(\Gamma)} S_{f} \otimes u_{s(e) s(f)} u_{r(e) r(f)}. $$ 
	Thus, by the universal property of $Q_{\tau}^{Lin}(\Gamma)$, there exists a CQG morphism $\phi: Q_{\tau}^{Lin}(\Gamma) \to QAut_{Ban}(\Gamma)$ sending $q_{ef} \mapsto u_{s(f) s(e)} u_{r(f) r(e)}$. Since the quantum symmetry of $C^*(\Gamma)$ is minimal, we have $q_{ee}^{*}q_{ee} = 1 = q_{ee}q_{ee}^{*}$. Therefore,
	\begin{align*}
	& \phi(q_{ee}^{*}q_{ee})=1=\phi(q_{ee}q_{ee}^{*}) \\
	\implies & u_{r(e) r(e)} u_{s(e) s(e)} u_{r(e) r(e)} = 1 = u_{s(e) s(e)} u_{r(e) r(e)} u_{s(e) s(e)}.
	\end{align*}
In other words, for each vertex $i,j \in V(\Gamma)$ with $i \to j$, we have:
 \begin{equation}
 u_{jj}u_{ii}u_{jj}=1   \label{Eq_ver1}
 \end{equation}
 \begin{equation}
 u_{ii}u_{jj}u_{ii}=1   \label{Eq_ver2}  
 \end{equation}
  Multiplying $u_{jj}$ to Equation \eqref{Eq_ver1} from the left, we get $u_{jj}=1$ (since $u_{jj}$ is a projection). Similarly, multiplying $u_{ii}$ to Equation \eqref{Eq_ver2} from the left, we get $u_{ii}=1$. \\ 
  Since $\Gamma$ has no isolated vertices, for every vertex $v \in V(\Gamma)$, at least one of $s^{-1}(v)$ or $r^{-1}(v)$ must be non-empty. Hence $u_{ii}=1$ for all $i \in V(\Gamma)$. Consequently, $\Gamma$ is always quantum asymmetric. \\
\end{proof}

\noindent To illustrate (2), consider the graph $P_{alt}$ shown in Figure \ref{P_alt}. Note that 
\begin{align*}
& indeg(v_1)=0,  \hspace{0.6cm} indeg(v_2)=2, \hspace{0.6cm} indeg(v_3)=0, \hspace{0.6cm} indeg(v_4)=1, \\
& outdeg(v_1)=1, \hspace{0.4cm} outdeg(v_2)=0, \hspace{0.4cm} outdeg(v_3)=2, \hspace{0.4cm} outdeg(v_4)=0.
\end{align*} 
 \begin{center}
	\begin{figure}[htb]
		\begin{tikzpicture}
		\draw[fill=black] (0,0) circle (2pt) node[anchor=south]{$v_{1}$};
		\draw[fill=black] (2,0) circle (2pt) node[anchor=south]{$v_{2}$};
		\draw[fill=black] (4,0) circle (2pt) node[anchor=south]{$v_{3}$};
		\draw[fill=black] (6,0) circle (2pt) node[anchor=south]{$v_{4}$};
		
		\begin{scope}[decoration={markings,
			mark=at position 0.50 with {\arrow[black,thick]{>}},
		}]
		\draw[black,postaction={decorate}, thick ](0,0)--node[above]{$e_{1}$}  (2,0);
		\end{scope}
		\begin{scope}[decoration={markings,
			mark=at position 0.50 with {\arrow[black,thick]{<}},
		}]
		\draw[black,postaction={decorate}, thick ](2,0)--node[above]{$e_{2}$}  (4,0);
		\end{scope}
		\begin{scope}[decoration={markings,
			mark=at position 0.50 with {\arrow[black,thick]{>}},
		}]
		\draw[black,postaction={decorate}, thick ](4,0)--node[above]{$e_{3}$}  (6,0);
		\end{scope}
		
		\end{tikzpicture}
		\caption{$P_{alt}$}   \label{P_alt}
	\end{figure} 
\end{center}   
By [Lemma 3.2.2, \cite{Fulton}], we know that if $indeg(v) \neq indeg(w)$ or $outdeg(v) \neq outdeg(w)$ for some $v,w \in V(\Gamma)$, then $u_{vw}=0$. Using this fact and comparing the indegree and outdegree of the vertices of $P_{alt}$, it is evident that $u_{v_{i}v_{j}}=0$ for all $i \neq j$. Consequently, $u_{v_{i}v_{i}}=1$ for all $i \in [4]$, and hence $P_{alt}$ is quantum asymmetric.\\[0.1cm]
To establish that $Q_{\tau}^{Lin}(P_{alt})$ is not minimal, we first note the Cuntz-Krieger relations for $C^*(P_{alt})$:
\begin{align*}
p_{v_1}= S_{e_1} S_{e_1}^*, \hspace{0.6cm} p_{v_2}=S_{e_1}^* S_{e_1} = S_{e_2}^* S_{e_2}, \hspace{0.6cm} p_{v_3}=S_{e_2}S_{e_2}^* + S_{e_3}S_{e_3}^*, \hspace{0.6cm} p_{v_4}=S_{e_3} S_{e_3}^*. 
\end{align*}
Observe that $\varphi:C^*(P_{alt}) \to C^*(P_{alt})$, defined on the generators by $\varphi(S_{e_1})=S_{e_2}$, $\varphi(S_{e_2})=S_{e_1}$ and $\varphi(S_{e_3})=S_{e_3}$, is an order 2 automorphism of $C^*(P_{alt})$.\\
We claim that $\left( \mathcal{D}_{\varphi}\left(\mystar\limits_{3} C(S^1)\right),\Delta_{\varphi}\right) $, the doubling of the CQG $\left(\mystar\limits_{3} C(S^1)\right)$ \cite{Soltan}, acts on $C^*(P_{alt})$ via a linear faithful and $\tau$-preserving action (say, $\alpha$). To that end, define an action $\alpha:C^*(P_{alt}) \to C^*(P_{alt}) \otimes \left( \mathcal{D}_{\varphi}\left(\mystar\limits_{3} C(S^1)\right),\Delta_{\varphi}\right) $  by 
\begin{align*}
\alpha(S_{e_1}) &= S_{e_1} \otimes (z_{1},0)+S_{e_2} \otimes (0,z_{2}),\\
\alpha(S_{e_2}) &= S_{e_1} \otimes (0,z_{1}) + S_{e_2} \otimes (z_{2},0),\\
\alpha(S_{e_3}) &= S_{e_3} \otimes (z_{3},z_{3}),
\end{align*}
where the unitary elements $\{z_i\}_{i=1}^{3}$ denote the canonical generators of $\mystar\limits_{3} C(S^1)$. Denote the action coefficient matrix $u=(u_{e_i e_j})_{3 \times 3}:=\begin{bmatrix}
    (z_1,0) & (0,z_1) & 0 \\
    (0,z_2) & (z_2,0) & 0 \\
    0 & 0 & (z_3,z_3)
\end{bmatrix}$. Now, extend $\alpha$ to the vertex projections $\{p_{i}\}_{i=1}^{4}$ using the Cuntz-Krieger relations as follows:
\begin{align*}
    \alpha(p_{v_1}) = & \alpha(S_{e_1} S_{e_1}^*) = p_1 \otimes (1,0) + S_{e_2}S_{e_2}^* \otimes (0,1), \\
    \alpha(p_{v_2}) = & \alpha(S_{e_1}^*S_{e_1}) = p_2 \otimes 1, \\
    \alpha(p_{v_3}) = & \alpha(S_{e_2}S_{e_2}^* + S_{e_3}S_{e_3}^*) = p_{v_1} \otimes (0,1) + S_{e_2} S_{e_2}^* \otimes (1,0) + S_{e_3} S_{e_3}^* \otimes (1,1), \\
    \alpha(p_{v_4}) = & \alpha(S_{e_3}^*S_{e_3}) = p_4 \otimes 1.
\end{align*}
Since $p_{v_2}= S_{e_1}^* S_{e_1} = S_{e_2}^* S_{e_2} $ and by the definition of $\alpha$, it is straightforward to verify that  $\alpha(S_{e_1}^* S_{e_1})=\alpha(S_{e_2}^* S_{e_2})$, the existence of $\alpha$ follows directly. Since $\alpha$ is linear and $u$ is unitary, the action equation, faithfulness and Poddles condition follow from the standard arguments of CGQ actions (see Lemmas 2.5.2 \& 2.5.3 in Chapter 2 of \cite{Ujjal_thesis}). Finally, $F^{P_{alt}}=I_{3}$ and $u$ is unitary imply $\alpha$ is $\tau$ preserving. \\ 

\begin{rem}
	In the above example, from the point of view of quantum symmetry of a graph, the trivial group is the only finite (classical) group that is a quantum subgroup of $QAut_{Ban}(P_{alt})$. More generally, $Aut(\Gamma)$ is the largest (classical) group that is a quantum subgroup of $QAut_{Ban}(\Gamma)$ for $\Gamma \in \mathcal{G}_{0,1}$ (see [Lemma 2.1.6, \cite{SSthesis}]). In contrast, from the graph $C^*$-algebraic point of view, $\mathbb{Z}_2$ (identifying its function algebra $C(\mathbb{Z}_2)$ with $\mathbb{C} \oplus \mathbb{C}$) acts on $C^*(P_{alt})$ via the linear faithful and $\tau$-preserving action defined by $S_{e_1} \mapsto S_{e_1} \otimes (1,0)+S_{e_2} \otimes (0,1)$, 
	$S_{e_2} \mapsto S_{e_1} \otimes (0,1) + S_{e_2} \otimes (1,0)$, and
	$S_{e_3} \mapsto S_{e_3} \otimes (1,1)$. Therefore, $\mathbb{Z}_2$ is a quantum subgroup of $Q_{\tau}^{Lin}(P_{alt})$. Thus, although $Aut(\Gamma)$ acts on $C^*(\Gamma)$ in general (by [Lemma 2.1.6, \cite{SSthesis}] and [Theorem 3.10, \cite{Mandal}]), it may not be the largest finite group that is contained in $Q_{\tau}^{Lin}(\Gamma)$ as a quantum subgroup. \\	
\end{rem}

Since $QAut_{Ban}(\Gamma)$ is a quantum subgroup of $Q_{\tau}^{Lin}(\Gamma)$, a natural question is whether the isomorphism of the quantum symmetry groups of two graphs implies the isomorphism of the quantum symmetry groups of their associated graph $C^*$-algebras, or vice versa. Clearly, there exist infinitely many pairs of finite weakly connected graphs in $\mathcal{G}_{0,1}$, say $\Gamma_1$ and $\Gamma_2$, with distinct numbers of edges having $QAut_{Ban}(\Gamma_1) \cong QAut_{Ban}(\Gamma_2)$ (by Theorem \ref{graph_infinitely many}). However, Corollary \ref{edge equal} ensures that $Q_{\tau}^{Lin}(\Gamma_1) \ncong Q_{\tau}^{Lin}(\Gamma_2)$.\\
The converse implication is also not true in general. For instance, take $\Gamma_1 = So_3$ (see Figure \ref{Son}) and $\Gamma_2 = P_1 \sqcup So_2$. Comparing the indegree and outdegree of the vertices and using [Lemma 3.2.2, \cite{Fulton}], it is straightforward to show that $QAut_{Ban}(So_3) \cong S_3^{+} \ncong S_2^+ \cong QAut_{Ban}(P_1 \sqcup So_2) $. In contrast, [Proposition 3.4, \cite{max_permutational_symm}] implies that $Q_{\tau}^{Lin}(So_3) \cong SH_3^{\infty +} \cong Q_{\tau}^{Lin}(P_1 \sqcup So_2)$.\\
However, in the family of graphs having minimal quantum symmetries of their associated graph $C^*$-algebras, invariance under $Q_{\tau}^{Lin}(.)$ implies invariance under $QAut_{Ban}(.)$. This is an immediate consequence of Proposition \ref{rigid imply q.assym} as shown in Corollary \ref{invarient under minimial}.
\begin{cor} \label{invarient under minimial}
    For $\Gamma_1, \Gamma_2 \in \mathcal{G}_{0,1}$, if $Q_{\tau}^{Lin}(\Gamma_1) \cong Q_{\tau}^{Lin}(\Gamma_2)$ and quantum symmetry of $C^*(\Gamma_1)$ is minimal, then $QAut_{Ban}(\Gamma_1) \cong QAut_{Ban}(\Gamma_2)$.
\end{cor}

\section*{Acknowledgements}
\noindent  The first author acknowledges the financial support from the Indian Statistical Institute Kolkata. The second author also acknowledges the support from DST-FIST (File No. SR/FST/MS-I/2019/41).

\end{document}